\documentclass[11pt,a4paper]{amsart}
\usepackage[T1]{fontenc}
\usepackage{lmodern}
\usepackage{amsmath,amssymb,amsfonts}
\usepackage[margin=28mm]{geometry}
\usepackage[hidelinks]{hyperref}
\DeclareMathOperator{\Hom}{Hom}
\DeclareMathOperator{\Spec}{Spec}
\DeclareMathOperator{\Supp}{Supp}
\DeclareMathOperator{\Fitt}{Fitt}
\DeclareMathOperator{\Lie}{Lie}
\DeclareMathOperator{\Ann}{Ann}
\DeclareMathOperator{\Ad}{Ad}
\newcommand{\OO}{\mathcal O}
\newcommand{\mcH}{\mathcal H}
\newcommand{\mcD}{\mathcal D}
\newcommand{\mcZ}{\mathcal Z}
\newcommand{\mcR}{\mathcal R}
\newcommand{\mcJ}{\mathcal J}
\newcommand{\mcI}{\mathcal I}
\newcommand{\sms}{\smallskip}
\newcommand{\Zspec}{\mathfrak Z^{\mathrm{spec}}}
\newcommand{\Qbar}{\overline{\mathbb Q}_\ell}

\newtheorem{theorem}{Theorem}[section]
\newtheorem{proposition}[theorem]{Proposition}
\newtheorem{corollary}[theorem]{Corollary}
\newtheorem{lemma}[theorem]{Lemma}
\theoremstyle{definition}
\newtheorem{definition}[theorem]{Definition}
\newtheorem{conjecture}[theorem]{Conjecture}
\theoremstyle{remark}
\newtheorem{remark}[theorem]{Remark}
\newtheorem*{remark*}{Remark}
\DeclareMathOperator{\GL}{GL}
\DeclareMathOperator{\tr}{tr}
\newcommand{\grs}{\mathfrak s}
\newcommand{\grt}{\mathfrak t}
\newcommand{\St}{\operatorname{St}}
\newcommand{\Icap}{\mathcal I^{\cap}}

\begin{document}
\title{Invariant distributions and common denominators\\
in the Bernstein center}
\author{Alexander Braverman and David Kazhdan}
\begin{abstract}
Let $G$ be a connected reductive group over a non-Archimedean local
field. We conjecture that every invariant distribution whose support is compact modulo conjugation belongs to the rational Bernstein center. We also make an explicit conjecture about the ideal of possible denominators (the description is using a version of the local Langlands correspondence and involves working with endoscopic subgroups of $G$). 

We prove this for products of general linear groups; in this case the above idea is much more explicit and Langlands correspondence is not needed for its description (it is a principle ideal, where the generator is given by some version of the discriminant function). We also prove our general conjecture for  $G=\mathrm{SL}(2)$.
\end{abstract}
\maketitle

\section{Introduction}
\subsection{Basic notations}
Let $F$ be a non-Archimedean local field of residue characteristic $p$
and residue-field cardinality $q$, and let $G/F$ be a connected
reductive group. Write $\OO_F$ for its ring of integers, fix a uniformizer
$\varpi$, and normalize $v_F$ and $|\cdot|$ by
$v_F(\varpi)=1$ and $|\varpi|=q^{-1}$.
Throughout, $G$ denotes the algebraic group and
$G(F)$ its locally compact group of $F$-points; we use the same
convention for $L$, $M$, $P$, and other algebraic groups. We denote by $\mathcal M(G)$ the category of smooth representations of $G(F)$
(in complex vector spaces).
\subsection{Invariant distributions and the Bernstein center}
Fix a Haar measure on $G(F)$. Write
$\mcH_G=C_c^\infty(G(F),\mathbb C)$ for its Hecke algebra and
$\mcD_G$ for the space of conjugation-invariant distributions on $G(F)$.
The subscripts in $\mcH_G$, $\mcD_G$, and $\mcZ_G$ refer to the
algebraic group $G$, while these spaces concern $G(F)$.
We realize the Bernstein center $\mcZ_G\subset\mcD_G$ as the algebra
of essentially compact invariant distributions. By definition these are invariant distributions $z$ such that  
$z*h\in\mcH_G$ for every $h\in\mcH_G$.

Let $\mcD_G^c$ consist of those invariant distributions $\alpha$ for
which there is a compact subset $C\subset G(F)$ such that
$$
 \Supp(\alpha)\subset C^{G(F)},\qquad
 C^{G(F)}=\{gcg^{-1}:g\in G(F),\ c\in C\}.
$$
Thus $\mcD_G^c$ is the space of invariant distributions with bounded
support modulo conjugation. The main examples of such distributions are {\em orbital measures} (a.k.a orbital integrals):
for every adjoint orbit $\Omega\subset G(F)$ it is known that 

a) $\Omega$ possesses a $G(F)$-ivariant measure $\mu$ (unique up to constant)

b) The integral $\int\limits_{\Omega} hd\mu$ is absolutely convergent. 

\noindent
Thus the above integral defines an invariant distribution which by definition lies in $\mcD^G_c$. 
\subsection{The rational Bernstein center}
The main goal of this paper 

Let
$\mcR_G$ be the rational Bernstein-center subspace:
$$
 \mcR_G=\{\alpha\in\mcD_G:
   z*\alpha\in\mcZ_G\text{ for some non-zero-divisor }z\in\mcZ_G\}.
$$
The central action on arbitrary distributions is the transpose of its
action on test functions. Explicitly, if $h^\vee(g)=h(g^{-1})$ and
$z^\vee(h)=z(h^\vee)$, then
$(z*\alpha)(h)=\alpha(z^\vee*h)$.

\begin{remark*}
Here rationality means saturation of $\mcZ_G$ inside $\mcD_G$ under
the action of non-zero-divisors. It does not assert that $\mcR_G$
embeds in the total fraction ring of $\mcZ_G$: an invariant
distribution annihilated by a non-zero-divisor already belongs to
$\mcR_G$. Such nonzero central-torsion distributions occur for
$\mathrm{SL}(2)$ in Section~\ref{sec:sl2}.
\end{remark*}

\begin{conjecture}\label{conj:denominator}
For every connected reductive group $G/F$ there exists a non-zero divisor $z\in \mcZ _G$ such that $z\ast \mcD_G^c \subset \mcZ _G$. \end{conjecture}
\begin{corollary}
$$
 {\mcD_G^c\subset\mcR_G.}
$$ \end{corollary}
\begin{remark}In other words, every invariant distribution  supported in the conjugates of a compact set belongs to the rational Bernstein center.\end{remark}

\subsection{The denominator for general linear groups}
We prove this conjecture for $G=\GL(N)$, and more generally for
products of general linear groups. More precisely, in this case we
can describe the entire ideal of common denominators
$$
 \mcJ_G=\{z\in\mcZ_G:z*\mcD_G^c\subset\mcZ_G\}.
$$
(The fact that this is an ideal is obvious). 

For $G=\GL(N)$, a Bernstein component $\grs$ is specified by a
cuspidal support
$$
 \left(\prod_a\GL(d_a,F)^{m_a},\ 
             \bigotimes_a\sigma_a^{\otimes m_a}\right),
 \qquad \sum_a d_am_a=N,
$$
where the $\sigma_a$ are pairwise inequivalent up to unramified twist.
Let $z_{a,i}$ be coordinates on their unramified-twist tori, after
quotienting by the finite groups of self-twists. Then
$$
 Z_{\grs}=\bigotimes_a
    \mathbb C[z_{a,1}^{\pm1},\ldots,z_{a,m_a}^{\pm1}]^{S_{m_a}},
 \qquad \mcZ_G=\prod_{\grs}Z_{\grs}.
$$
Define the discriminant by
\begin{equation}\label{eq:intro-discriminant}
 \Delta_{\grs}
   =\prod_a\prod_{\substack{1\leq i,j\leq m_a\\i\ne j}}
            \left(1-\frac{z_{a,i}}{z_{a,j}}\right),
 \qquad \Delta_G=(\Delta_{\grs})_{\grs}\in\mcZ_G.
\end{equation}
The empty product is $1$. The definition is independent of the chosen
cuspidal representatives. For a product $L$ of general linear
groups, define $\Delta_L$ by taking the product of the discriminants
for its factors.

\begin{theorem}\label{T:denominator-exact}
For every product $L$ of general linear groups,
$$
 {\mcJ_L=\Delta_L\mcZ_L.}
$$
The element $\Delta_L$ is a non-zero-divisor. In particular,
$$
 \Delta_L*\mcD_L^c\subset\mcZ_L,
 \qquad \mcD_L^c\subset\mcR_L.
$$
Thus Conjecture~\ref{conj:denominator} holds for $L$, and every
common denominator is divisible by $\Delta_L$ component by component.
\end{theorem}

\subsection{The endoscopic denominator ideal}
For general $G$, our proposed description of $\mcJ_G$ uses
Langlands parameters. We briefly explain the objects involved; precise
algebraic constructions are given in Sections~\ref{sec:spectral}--\ref{sec:endoscopy}.

Let $\widehat G$ be the complex Langlands dual group of $G$, and let
$W_F$ be the Weil group of $F$. The $L$-group is
${}^LG=\widehat G\rtimes W_F$, where $W_F$ acts on $\widehat G$
through the action determined by $G$. A (Weil) Langlands parameter
for $G$ is a continuous homomorphism
$$
 \varphi:W_F\longrightarrow{}^LG
$$
whose composition with the projection ${}^LG\to W_F$ is the identity.
Two parameters are identified if they are conjugate by $\widehat G$.
When $G$ is split, the Weil action is trivial, so a parameter is simply
a continuous homomorphism $W_F\to\widehat G$, up to conjugation.
We use Weil parameters here, rather than Weil--Deligne parameters;
no additional monodromy operator is included.

The ring $A_G=\Zspec_G$, called the \emph{spectral Bernstein center},
consists of regular functions on the stack of these parameters.
Fargues--Scholze construct a homomorphism
$$
 \Psi_G:A_G\longrightarrow\mcZ_G.
$$
The parameter space and this map are defined in
Section~\ref{sec:spectral}.

An ordinary untwisted endoscopic datum $\mathfrak e$ for $G$ includes
another quasi-split reductive group $H$ and a group
$H'$ fitting into an exact sequence
$$
 1\longrightarrow\widehat H\longrightarrow H'
 \xrightarrow{p_{\mathfrak e}}W_F\longrightarrow1,
$$
together with an embedding
$$
 \eta:H'\longrightarrow{}^LG
$$
compatible with the projections to $W_F$. In particular, its restriction
to $\widehat H$ identifies $\widehat H$ with a subgroup of $\widehat G$.
When an appropriate $L$-embedding ${}^LH\to{}^LG$ has been chosen,
we can take $H'={}^LH$.
A parameter with values in $H'$ is a continuous
homomorphism $\varphi_H:W_F\to H'$ satisfying
$p_{\mathfrak e}\circ\varphi_H=\operatorname{id}_{W_F}$.
Composing $\varphi_H$ with $\eta$ gives a parameter for $G$.
This is the \emph{transfer of parameters} attached to $\mathfrak e$.

For such a parameter, let $Z_{\widehat H}(\varphi_H)$ denote the
subgroup of $\widehat H$ commuting with the image of the parameter
(in $H'$), and similarly for
$Z_{\widehat G}(\eta\circ\varphi_H)$. The first centralizer is
contained in the second. We call $\varphi_H$ \emph{exceptional} if
$$
 \dim Z_{\widehat G}(\eta\circ\varphi_H)
       >\dim Z_{\widehat H}(\varphi_H).
$$
Thus a parameter is exceptional precisely when its transferred
parameter has additional infinitesimal symmetries not present for $H$.

To turn this set into an ideal, we retain more information than just
whether the inequality holds. Put
$V=\Lie(\widehat G)/\Lie(\widehat H)$. The parameter acts linearly
on $V$ through $\eta$ and the adjoint action. On a finite-type
parameter space, this action factors through a finitely generated
group $\Gamma$; choose generators $\gamma_1,\ldots,\gamma_m$.
For each parameter form the linear map
$$
 T_{\varphi_H}:V\longrightarrow V^{\oplus m},\qquad
 v\longmapsto\bigl((\gamma_1-1)v,\ldots,(\gamma_m-1)v\bigr).
$$
The exceptional condition is equivalent to $\ker T_{\varphi_H}\ne0$.
Writing $n=\dim V$, take the ideal generated by all $n\times n$
minors of the matrix of $T$ (and take the unit ideal if $n=0$).
Its vanishing set is exactly the exceptional locus. The minors also
record the order of vanishing along that locus; consequently we do
not replace this ideal by the ideal of functions merely vanishing on
the underlying set. The construction does not depend on the chosen
generators or bases; see Section~\ref{sec:algebraic}.

We pass from the parameter space to its quotient by $\widehat H$:
intersect the ideal of minors with the ring of invariant regular
functions. On each bounded parameter piece, denote the resulting
ideal by $J_{\mathfrak e,U}$, and the induced transfer map by
$$
 f_{\mathfrak e,U}:Y_{\mathfrak e,U}\longrightarrow Y_{G,U}.
$$
The condition imposed on a regular function $a$ on the target is
simply that its pullback belongs to $J_{\mathfrak e,U}$.
In symbols, the corresponding ideal on the target is
$$
 \mcI_{G,U}^{\mathfrak e}
   =\{a\in A_{G,U}:f_{\mathfrak e,U}^*a\in J_{\mathfrak e,U}\}.
$$
Intersect these ideals over all ordinary untwisted endoscopic data,
and impose the same condition on every bounded piece. This defines
$\Icap_G\subset A_G$; see \eqref{eq:global-intersection}.
The definition is entirely in terms of regular functions and explicit
matrix minors, although the more invariant formulation used later is
convenient for proofs.

\begin{conjecture}[Endoscopic description of the denominator ideal]
\label{conj:geometric}
The endoscopic ideal $\Icap_G$ is intrinsic and satisfies
\begin{equation}\label{eq:intersection-equality}
 {\mcJ_G=\mcZ_G\cdot\Psi_G(\Icap_G).}
\end{equation}
This ideal contains a non-zero-divisor of $\mcZ_G$.
\end{conjecture}
This refines Conjecture~\ref{conj:denominator} by describing the ideal
of common denominators.

We also prove Conjecture~\ref{conj:geometric} for $\GL(N)$. Under
the spectral-to-Bernstein identification, the component of
$\Icap_G$ on $\grs$ is exactly $(\Delta_{\grs})$; see
Theorem~\ref{T:intersection-typeA}. Thus the endoscopic construction
recovers the denominator ideal of Theorem~\ref{T:denominator-exact}.

\subsection{Sketch of the proof of Theorem~\ref{T:denominator-exact}}
Roughly speaking, the idea of the proof is this. Fix a congruence subgroup $K$ of $G(F)$. Then one can show that the restriction of any invariant distribution with support bounded modulo conjugation to the corresponding Hecke algebra is a sum of distributions induced from Bernstein center elements of Levi subgroups. Comparing type $A$
Plancherel measures shows that the discriminant kills their
denominators. For the reverse
inclusion, an induced cuspidal Bernstein projector has a
numerator relatively prime to the discriminant.

\subsection{Outline of the paper}
We first give the precise general construction of the endoscopic
ideal in Sections~\ref{sec:spectral}--\ref{sec:endoscopy}.
Section~\ref{sec:typeA} proves Theorem~\ref{T:denominator-exact}.
Section~\ref{sec:intersection-typeA} computes the endoscopic
intersection.
Section~\ref{sec:sl2} proves both conjectures for $\mathrm{SL}(2)$.
In that case the denominator is $t^2-4$ on each principal series
component with inversion symmetry, where $t=z+z^{-1}$, and is $1$
on the other components. 

\subsection{Acknowledgements} We thank R.~Bezrukavnikov and M.~Finkelberg for helpful discussions. During the preparation of this paper, D.K. was partially supported by ERC grant no.~101142781. A.B. was partially supported by NSERC. This work was done while A.B. was visiting the Simons Laufer Mathematical Institute. 

\section{The spectral Bernstein center}\label{sec:spectral}
Fix a prime $\ell\ne p$ and put $\Lambda=\Qbar$. Let $\widehat G$ be the Langlands dual group over $\Lambda$, equipped with the pinned action of $W_F$. Fix an isomorphism $\iota:\Lambda\simeq\mathbb C$ and a square root of $q$ used in the Fargues--Scholze normalization. All comparisons with complex distributions use these choices.
\subsection{The parameter scheme}
Fargues--Scholze construct a scheme $Z^1(W_F,\widehat G)_\Lambda$ representing condensed $1$-cocycles; see \cite[Definition VIII.1.1 and Theorem VIII.1.3]{FS}. A cocycle $c$ satisfies
$$
c(ww')=c(w)\,w(c(w')),\qquad w,w'\in W_F,
$$
together with the condensed continuity condition. Equivalently, it defines a section
$$
\varphi_c:W_F\longrightarrow\widehat G\rtimes W_F,\qquad w\longmapsto(c(w),w)
$$
of the projection to $W_F$. For split $G$, the Weil action is trivial and these are continuous homomorphisms into $\widehat G$.
The conjugation action is
$$
(g\cdot c)(w)=g\,c(w)\,w(g)^{-1}.
$$
The parameter stack and spectral Bernstein center are
$$
\operatorname{Loc}_G=[Z^1(W_F,\widehat G)_\Lambda/\widehat G],
\qquad
A_G:=\Zspec_G=\Gamma(\operatorname{Loc}_G,\OO).
$$
Thus
$$
A_G=\OO(Z^1(W_F,\widehat G)_\Lambda)^{\widehat G}.
$$
The full parameter space is handled using finite-type pieces with bounded wild ramification. On characteristic-zero affine quotients, closed geometric points correspond to semisimple parameter classes. In particular, for an invariant regular function $a$,
$$
a(\varphi)=a(\varphi^{\mathrm{ss}}).
$$
\subsection{The map to the ordinary Bernstein center}
\begin{theorem}[Fargues--Scholze]
With the coefficient and normalization choices above, there is a ring homomorphism
$$
\Psi_G:A_G\longrightarrow\mcZ_G.
$$
\end{theorem}
The construction is given in \cite[Proposition I.9.3, Theorem IX.5.2, and Definition IX.7.1]{FS}.
\subsection{Stable and very stable centers}\label{sec:stable-centers}
We spell out the harmonic-analytic definitions used in
\cite[Section 1.1]{H} and \cite[Definition 1.1.1]{Varma}.
Two strongly regular semisimple elements of $G(F)$ are stably
conjugate if they are conjugate under $G(\overline F)$.
For such an element $\gamma\in G(F)$, let $G_\gamma$ denote its
algebraic centralizer, which is a torus. For $f\in\mcH_G$, put
$$
 O_\gamma(f)=\int_{G_\gamma(F)\backslash G(F)}
                f(x^{-1}\gamma x)\,d\dot x.
$$
Choose Haar measures on the centralizer tori compatibly under stable
conjugacy. The stable orbital integral is
$$
 SO_\gamma(f)=\sum_{[\gamma']\subset[\gamma]_{\mathrm{st}}}
                         O_{\gamma'}(f),
$$
where the finite sum runs over the $G(F)$-conjugacy classes in the
stable conjugacy class of $\gamma$ that occur in $G(F)$.
Changing the compatible measures multiplies this functional by a
nonzero scalar and does not change its kernel.

Let
$$
 \mcH_G^{\mathrm{unst}}
  =\{f\in\mcH_G:SO_\gamma(f)=0
       \text{ for every strongly regular semisimple }\gamma\in G(F)\}.
$$
An invariant distribution $D$ on $G(F)$ is \emph{stable} if
$D(f)=0$ for every $f\in\mcH_G^{\mathrm{unst}}$.
In particular, the \emph{stable Bernstein center} is the vector space
$$
 \mcZ_G^{\mathrm{st}}
   =\{z\in\mcZ_G:z(f)=0\text{ for all }f\in\mcH_G^{\mathrm{unst}}\}.
$$
Here stability is a condition on $z$ viewed as a distribution;
closure of this space under convolution is not part of the definition.

The \emph{very stable Bernstein center} is
$$
 \mcZ_G^{\mathrm{vst}}
   =\{z\in\mcZ_G:z*\mcH_G^{\mathrm{unst}}
                              \subset\mcH_G^{\mathrm{unst}}\}.
$$
Thus $z$ is very stable precisely when convolution by $z$ descends
to an operator on the stable quotient
$\mcH_G/\mcH_G^{\mathrm{unst}}$. Unlike the definition of stability,
this condition concerns the action on all unstable test functions.
It immediately makes $\mcZ_G^{\mathrm{vst}}$ a subalgebra of
$\mcZ_G$.

For $p$-adic $F$, one has
$\mcZ_G^{\mathrm{vst}}\subset\mcZ_G^{\mathrm{st}}$;
see \cite[Remark 1.1.2]{Varma}. One way to see the implication uses
the stability of the distribution $f\mapsto f(1)$, together with
$z(f)=(z*f^\vee)(1)$ and the preservation of
$\mcH_G^{\mathrm{unst}}$ by inversion. Moreover,
$\mcZ_G^{\mathrm{st}}$ is a module over
$\mcZ_G^{\mathrm{vst}}$ under convolution.

\begin{remark}
For $p$-adic $F$, the equality
$$
 \mcZ_G^{\mathrm{st}}=\mcZ_G^{\mathrm{vst}}
$$
is conjectured in general; see \cite[Conjecture 1.1.4]{Varma}
and \cite[Section 1.1]{H}. It is stronger than merely asserting
that stable central distributions form a subalgebra. Varma proves
the equality under the hypotheses on tempered $L$-packets specified
in \cite[Hypothesis 2.5.1 and Theorem 4.4.2]{Varma}.
Hansen's theorem below gives an unconditional inclusion for the
image of the spectral center; it does not assert equality of the
stable and very stable centers for arbitrary $G$.
\end{remark}

\begin{theorem}[Hansen]
If $F$ is a finite extension of $\mathbb Q_p$, then
$$
\operatorname{Im}(\Psi_G)\subset\mcZ_G^{\mathrm{vst}}\subset\mcZ_G^{\mathrm{st}}.
$$
\end{theorem}
The assertion is \cite[Theorem 1.1]{H}.
\begin{remark}
The cited theorem of Hansen assumes characteristic zero. No stability assertion for fields of positive characteristic is deduced from that reference here. The definitions and conjectures in this article are nevertheless formulated for arbitrary non-Archimedean local fields.
\end{remark}
\section{A general algebraic construction}\label{sec:algebraic}
We first give a construction for ordinary representation varieties. Its adaptation to endoscopic parameter spaces is described in Section~\ref{sec:endoscopy}.

\sms

\subsection{Setup}
Let $k$ be a field of characteristic zero, let $\hat G$ be a reductive algebraic
group over $k$, and let $i: \hat H\hookrightarrow \hat G$ be a closed reductive subgroup.
Let $\Gamma$ be a finitely generated discrete group. Write
$$
 A(K)=\Hom(\Gamma,K),\qquad
 B(K)=A(K)//K=\Spec\bigl(k[A(K)]^K\bigr),\quad K= \hat H, \hat G.
$$
We use reduced representation varieties; the same arguments also apply to
the affine representation schemes. The induced morphism
$$
 f:B(\hat H)\longrightarrow B(\hat G)
$$
is finite by the Vinberg finiteness theorem (see \cite[Theorem 1.1]{Cotner}).
In particular, its image $C$ is closed.

For a geometric point $c'\in B(\hat H)$, choose a representation
$\rho:\Gamma\to \hat H$ in the unique closed $\hat H$-orbit over $c'$, and put
$$
 d(c')=\dim Z_{\hat G}(i\circ\rho)-\dim Z_{ \hat H}(\rho).
$$
Here $Z_K(\rho)$ denotes the centralizer of the image of the representation,
equivalently its stabilizer under conjugation. This integer is independent
of the choice within that closed orbit. For $r\geq1$, define
$$
 Z_r=\left\{c\in C:\max_{c'\in f^{-1}(c)}d(c')\geq r\right\},
 \qquad Z=Z_1.
$$
Fibers and maxima are interpreted over an algebraic closure of the residue
field. The fibers are finite and nonempty for $c\in C$.

\subsection{Closedness of the exceptional locus}
\begin{proposition}
For every integer $r\geq1$, the subset $Z_r$ is Zariski closed in $B(\hat G)$.
\end{proposition}
\begin{proof}
We first work over an algebraically closed field. Set
$\mathfrak g=\Lie(\hat G)$ and $\mathfrak h=\Lie(\hat H)$.
Since $\hat H$ is linearly reductive in characteristic zero, the exact sequence
of $\hat H$-modules
$$
0\longrightarrow\mathfrak h\longrightarrow\mathfrak g
\longrightarrow V\longrightarrow0,\qquad V=\mathfrak g/\mathfrak h,
$$
splits $\hat H$-equivariantly. Centralizers are smooth in characteristic zero,
and their Lie algebras are the corresponding fixed subspaces. Hence,
for every $\rho\in A(\hat H)$, including representations with nonclosed orbits,
\begin{equation}\label{eq:fixed}
 \dim Z_{\hat G}(i\circ\rho)-\dim Z_{\hat H}(\rho)
 =\dim\mathfrak g^{\rho(\Gamma)}-\dim\mathfrak h^{\rho(\Gamma)}
 =\dim V^{\rho(\Gamma)}.
\end{equation}

Choose generators $\gamma_1,\ldots,\gamma_m$ of $\Gamma$. The fixed
subspace in \eqref{eq:fixed} is the kernel of
$$
 T_\rho:V\longrightarrow V^{\oplus m},\qquad
 v\longmapsto\bigl((\rho(\gamma_j)-1)v\bigr)_{j=1}^m,
$$
where $\hat H$ acts on $V$ through the adjoint action. Its matrix entries depend
regularly on $\rho$. Thus
$$
 D_r=\{\rho\in A(\hat H):\dim\ker T_\rho\geq r\}
$$
is a closed, $\hat H$-invariant determinantal subset.

Let $q_{\hat H}:A( \hat H)\to B(\hat H)$ be the quotient map. An affine reductive quotient
sends invariant closed subsets to closed subsets; see
\cite[Lemma 2.3]{BGM}. Moreover,
\begin{equation}\label{eq:quotientlocus}
 q_{\hat H }(D_r)=\{c'\in B(\hat H):d(c')\geq r\}.
\end{equation}
Indeed, if a quotient fiber meets $D_r$, then the closure of the orbit of
such a point is contained in $D_r$. That orbit closure contains the unique
closed orbit of the fiber. Conversely, a closed-orbit representative with
$d(c')\geq r$ belongs to $D_r$ by \eqref{eq:fixed}.

Since $f$ is finite and hence closed, \eqref{eq:quotientlocus} gives
$$
 Z_r=f\bigl(q_{\hat H}(D_r)\bigr),
$$
which is closed in $B(\hat G)$. All constructions are defined over $k$;
invariants commute with field extension in characteristic zero. Thus the
argument gives the stated closed locus over the original field as well.
\end{proof}

\subsection{A canonical coherent sheaf}
The following construction provides a natural
scheme structure obtained from a fixed-vector rank condition.

Put
$$
 R=k[A(\hat H)],\qquad S=R^{\hat H},\qquad V=\mathfrak g/\mathfrak h.
$$
The universal representation $\rho_{\mathrm{univ}}:\Gamma\to \hat H(R)$
acts on $N=R\otimes_k V^\vee$. Form its module of coinvariants
\begin{equation}\label{eq:coinvariants}
 M=N_\Gamma
 =N\big/\langle\gamma u-u:\gamma\in\Gamma,\ u\in N\rangle_R.
\end{equation}
This definition is canonical: neither generators of $\Gamma$ nor a
complement to $\mathfrak h$ in $\mathfrak g$ are part of the construction.
The module is finitely presented. Indeed, for any finite generating set
$\gamma_1,\ldots,\gamma_m$, it is the cokernel of
$$
 N^{\oplus m}\longrightarrow N,\qquad
 (u_1,\ldots,u_m)\longmapsto
 \sum_{j=1}^m(\gamma_j-1)u_j.
$$
At a geometric point $\rho$, coinvariants commute with specialization and
finite-dimensional duality gives
\begin{equation}\label{eq:fiber}
 M\otimes_R k(\rho)
 \simeq(V^\vee\otimes_k k(\rho))_{\rho(\Gamma)}
 \simeq\bigl((V\otimes_k k(\rho))^{\rho(\Gamma)}\bigr)^\vee.
\end{equation}

Choose generators $\gamma_1,\ldots,\gamma_m$ of $\Gamma$ and a
basis of $V$. The relations in \eqref{eq:coinvariants} are represented
by an $n\times mn$ matrix, where $n=\dim V$; its blocks are the
matrices of $\gamma_j-1$ acting on $V^\vee$.
Let $I$ be the ideal generated by all its $n\times n$ minors (with
$I=R$ if $n=0$). Although this description uses choices, the ideal
is independent of them: it is the standard zeroth determinantal
(or Fitting) ideal of the module $M$. Define
$$
 I=\Fitt_0^R(M),\qquad J=I\cap S,
 \qquad E=\Spec(S/J)\subset B(\hat H).
$$
The universal construction is equivariant for conjugation by $\hat H$,
so $I$ is $\hat H$-stable. The zero set of $I$ is $\Supp_R(M)$, which by
\eqref{eq:fiber} is precisely $D_1$.
Exactness of invariants in characteristic zero gives
$$
 S/J\simeq(R/I)^{\hat H}.
$$
Consequently $E$ is the affine quotient of the canonical determinantal
subscheme $\Spec(R/I)$, and its underlying set is
$$
 |E|=q_{\hat H}(D_1)=\{c'\in B(\hat H):d(c')>0\}.
$$

\begin{proposition}
Regard $\mathcal O_E$ as a coherent sheaf on $B(\hat H)$ via its closed
inclusion. Then
$$
 {\mathcal F=f_*\mathcal O_E}
$$
is a canonically defined coherent sheaf on $B(\hat G)$, and
$\Supp(\mathcal F)=Z$.
\end{proposition}
\begin{proof}
The sheaf $\mathcal O_E$ is coherent, and pushforward along the finite
morphism $f$ preserves coherence. For a finite morphism, the support of
the pushforward of a coherent sheaf is the image of its support. Thus
$$
 \Supp(\mathcal F)=f(|E|)=Z.
$$
\end{proof}

More explicitly, if $T=k[B( \hat G)]$ and $f^*:T\to S$ is the restriction map,
then $\mathcal F$ corresponds to the finite $T$-module $S/J$. Its annihilator
is
$$
 \operatorname{Ann}_T(S/J)=(f^*)^{-1}(J).
$$
Hence this construction also defines a canonical closed subscheme of
$B( \hat G)$ with underlying set $Z$.

\begin{remark}
The complement used in the proof of \eqref{eq:fixed} is not canonical,
but the module $V=\mathfrak g/\mathfrak h$, the coinvariant module $M$,
its Fitting ideal, and the sheaf $\mathcal F$ are canonical.
Replacing $E$ by $E_{\mathrm{red}}$ gives another canonical coherent sheaf
with the same support.
\end{remark}

\section{Endoscopic parameter spaces and the candidate ideal}\label{sec:endoscopy}
\subsection{A discrete Weil-group replacement}
Write $P_F\subset W_F$ for wild inertia. Let $U\subset P_F$ be open in $P_F$ and normal in $W_F$. The quotient $W_F/U$ still contains profinite tame inertia, so it is not a finitely generated discrete group.
Choose an arithmetic Frobenius $\mathrm{Fr}$ and a pro-generator $s_t$ of tame inertia $I_F/P_F$. Let $W_F^0$ be the inverse image in $W_F$ of
$$
\langle\mathrm{Fr},s_t\rangle
\simeq\mathbb Z[1/q]\rtimes\mathbb Z\subset W_F/P_F,
\qquad \mathrm{Fr}s_t\mathrm{Fr}^{-1}=s_t^q.
$$
Endow $W_F^0$ with the topology for which $P_F$ is open with its profinite topology and $W_F^0/P_F$ is discrete. Then
$$
\Gamma_U=W_F^0/U
$$
is a finitely generated discrete group: it is an extension of the two-generator group $\langle\mathrm{Fr},s_t\rangle$ by the finite group $P_F/U$.
This is the discretized Weil group used in \cite[Introduction and Sections 2--4]{DHKM}. For sufficiently small $U$, the Weil action on the dual group factors through $\Gamma_U$, and the bounded parameter scheme is described by
$$
X_{G,U}=Z^1(\Gamma_U,\widehat G),\qquad
Y_{G,U}=X_{G,U}//\widehat G,\qquad
A_{G,U}=\OO(Y_{G,U}).
$$
These are cocycle schemes, not ordinary homomorphism schemes unless the action is trivial. The characteristic-zero comparison with the condensed parameter construction is used to regard these quotients as bounded pieces of the spectral parameter space; see \cite[Section VIII.1]{FS} and \cite{DHKM}.
\subsection{Endoscopic data}
Fix an endoscopic datum
$\mathfrak e=(H,H',s_e,\eta)$, where $H$ is
quasi-split and $H'$ fits into an exact sequence
$$
1\longrightarrow\widehat H\longrightarrow H'
\xrightarrow{p_{\mathfrak e}}W_F\longrightarrow1.
$$
The embedding
$$
\eta:H'\longrightarrow{}^LG,
\qquad
\eta(\widehat H)=Z_{\widehat G}(s_e)^\circ,
$$
is compatible with the projections to $W_F$.
Here $s_e\in\widehat G$ is semisimple, and elements of $\eta(H')$ conjugate $s_e$ to itself modulo $Z(\widehat G)$. When an $L$-embedding ${}^LH\to{}^LG$ has been chosen, one can use $H'={}^LH$. In general the extension $H'$ must be retained; a bare inclusion $\widehat H\subset\widehat G$ does not specify parameter transfer.
Let $\varphi_H:W_F\to H'$ be a continuous homomorphic section of $p_{\mathfrak e}$. Define
$$
d_{\mathfrak e}(\varphi_H)
=\dim Z_{\widehat G}(\eta\circ\varphi_H)
 -\dim Z_{\widehat H}(\varphi_H).
$$
The centralizers use the full $L$-parameters, including their Weil action.
Put $\mathfrak g=\Lie(\widehat G)$, $\mathfrak h=\Lie(\widehat H)$, and $V_{\mathfrak e}=\mathfrak g/\mathfrak h$. Since
$$
\mathfrak h=\mathfrak g^{\Ad(s_e)},
$$
the eigenspace decomposition of $\Ad(s_e)$ gives a complement to $\mathfrak h$ preserved by $\eta(H')$. Indeed, conjugating $s_e$ by such an element changes it only by a central factor, which acts trivially on $\mathfrak g$. Consequently the splitting is compatible with every $\varphi_H$, and
\begin{equation}\label{eq:endofixed}
d_{\mathfrak e}(\varphi_H)
=\dim V_{\mathfrak e}^{\Ad(\eta\circ\varphi_H)(W_F)}.
\end{equation}
This supplies the analogue of \eqref{eq:fixed} in the endoscopic setting.
\subsection{The ideal on a bounded parameter piece}
Choose compatible wild-ramification bounds for the source and target, so that parameter transfer gives an algebraic morphism
$$
f_{\mathfrak e,U}:Y_{\mathfrak e,U}\longrightarrow Y_{G,U}.
$$
Here $X_{\mathfrak e,U}$ is the affine scheme of sections into the endoscopic extension with the chosen bounds, and $Y_{\mathfrak e,U}=X_{\mathfrak e,U}//\widehat H$. The choices must ensure that the universal adjoint action on $V_{\mathfrak e}$ factors through a finitely generated discrete Weil-group quotient. Write this quotient as $\Gamma$.
\begin{remark}\label{rem:finite}
The coherent-pushforward construction below requires $f_{\mathfrak e,U}$ to be finite. Section~\ref{sec:algebraic} proves this in the ordinary representation-variety situation by Cotner's theorem. Its use for an endoscopic cocycle model requires the corresponding finiteness result for that model. This is an explicit requirement of the endoscopic construction, rather than a consequence asserted solely from $\widehat H\subset\widehat G$.
\end{remark}
Assume this finiteness condition, and put
$$
R_{\mathfrak e,U}=\OO(X_{\mathfrak e,U}),\qquad
S_{\mathfrak e,U}=R_{\mathfrak e,U}^{\widehat H}.
$$
Using the universal adjoint action, define
$$
M_{\mathfrak e,U}
=(R_{\mathfrak e,U}\otimes_\Lambda V_{\mathfrak e}^\vee)_\Gamma,
\quad
K_{\mathfrak e,U}=\Fitt_0^{R_{\mathfrak e,U}}(M_{\mathfrak e,U}),
\quad
J_{\mathfrak e,U}=K_{\mathfrak e,U}\cap S_{\mathfrak e,U}.
$$
Set
$$
E_{\mathfrak e,U}=\Spec(S_{\mathfrak e,U}/J_{\mathfrak e,U}),
\qquad
\mathcal F_{\mathfrak e,U}=(f_{\mathfrak e,U})_*\OO_{E_{\mathfrak e,U}}.
$$
The splitting in \eqref{eq:endofixed} and the same determinantal and quotient arguments as in Section~\ref{sec:algebraic} show that $E_{\mathfrak e,U}$ has underlying set
$$
\{[\varphi_H]:d_{\mathfrak e}(\varphi_H)>0\},
$$
where closed-orbit representatives are used. Define the spectral ideal
\begin{equation}\label{eq:endoideal}
{\mcI_{G,U}^{\mathfrak e}
=\Ann_{A_{G,U}}(\mathcal F_{\mathfrak e,U})
=(f_{\mathfrak e,U}^*)^{-1}(J_{\mathfrak e,U})\subset A_{G,U}.}
\end{equation}
Its zero set is the image of the exceptional locus. Formula~\eqref{eq:endoideal} specifies the scheme structure, not just its underlying set.
\begin{remark}
Replacing this annihilator by its radical, or by the zeroth Fitting ideal of the pushforward, can change the ideal and hence the conjecture. We use exactly the annihilator in \eqref{eq:endoideal}.
\end{remark}
\subsection{The intersection ideal}
Use ordinary untwisted endoscopic data, with their realized
parameter-transfer maps, taken up to isomorphism and conjugacy over
the target. The usual untwisted endoscopic condition is part of this
convention; transfers with a nontrivial twisting character are not
included. For a collection $\mathcal E_U$ on a bounded piece, put
\begin{equation}\label{eq:intersection-ideal}
 {\Icap_{G,U}
   =\bigcap_{\mathfrak e\in\mathcal E_U}\mcI_{G,U}^{\mathfrak e}.}
\end{equation}
The empty intersection is the unit ideal. For a finite collection,
this can also be written
$$
 \Icap_{G,U}
 =\Ann_{A_{G,U}}\left(
       \bigoplus_{\mathfrak e\in\mathcal E_U}
                  \mathcal F_{\mathfrak e,U}\right).
$$
In particular, repeating an equivalent transfer does not change the
ideal. Distinct transfers with the same annihilator also impose just
one ideal-theoretic condition.

Here is a formulation that does not assume that an arbitrary
intersection commutes with localization.
On every bounded target piece consider all compatible bounded source
pieces and all endoscopic transfers as above. Define
\begin{equation}\label{eq:global-intersection}
 \Icap_G=\left\{a\in A_G:
 f_{\mathfrak e,U}^*(a|_{Y_{G,U}})\in J_{\mathfrak e,U}
 \text{ for every such }(\mathfrak e,U)\right\}.
\end{equation}
This is an ideal, defined by simultaneous pullback conditions.
Equivalently, its restrictions satisfy all the annihilator conditions
in \eqref{eq:endoideal}. The pullback definition makes sense even
before finiteness of the transfers is established; when they are
finite it has the coherent-sheaf interpretation given above.
For a finite collection on a finite-type piece, the local construction
is exactly \eqref{eq:intersection-ideal}. Compatibility and independence
of the bounded parameter models are included in the geometric
Conjecture~\ref{conj:geometric}. They hold in the general-linear computation in
Section~\ref{sec:intersection-typeA}.

This completes the definition of the ideal appearing in
Conjecture~\ref{conj:geometric}. Its proposed image in the Bernstein
center is the common-denominator ideal $\mcJ_G$.

\section{Proof of the theorem for general linear groups}\label{sec:typeA}

\subsection{Some notation} In this section $G=\GL(N)$ and $L=\prod_{j=1}^s\GL(n_j)$.
We fix Haar measures on $L(F)$ and on $M(F)$ for every Levi subgroup $M$ of $L$.
Write
$$
 \mcH_L=C_c^\infty(L(F)),\qquad
 \overline\mcH_L=
 \mcH_L/[\mcH_L,\mcH_L],
$$
and let $\mcZ_L$ be the Bernstein center. For a compact open
subgroup $K\subset L(F)$, put
$$
 e_K=\operatorname{vol}(K)^{-1}\mathbf1_K,\qquad
 \mcH(L(F),K)=e_K*\mcH_L* e_K.
$$
An invariant distribution is a functional on $\overline\mcH_L$.
We regard every $b\in\mcZ_L$ as an invariant distribution and write
$b(h)$ for its value on a test function, or on its class in the cocenter.
We also use the realization of $b$ as a regular function on the Bernstein
variety. These descriptions are related by the Plancherel formula
$$
 b(h)=\int_{\widehat{L(F)}_{\mathrm{temp}}}
           b(\pi)\,\tr\pi(h^\vee)\,d\nu_L(\pi).
$$
Here $\nu_L$ is the Plancherel measure. Essential compactness means that $b*e_K\in\mcH(L(F),K)$ for every $K$.

\subsection{Bounded functionals and maps attached to Levi subgroups}\label{SS:bern-auxiliary}

We say that a subset of $L(F)$ is \emph{bounded modulo conjugation} if
it is contained in
$$
 L(F)\cdot C=\{gcg^{-1}:g\in L(F),\ c\in C\}
$$
for some compact subset $C\subset L(F)$. 

Here is a convenient family of sets exhausting all such supports.
Extend the absolute value of $F$ to an algebraic closure. For an integer
$R\geq0$, let
\begin{equation}\label{E:bern-bounded-set}
 B_R(L(F))=\{g=(g_j)\in L(F):
     q^{-R}\leq|\alpha|\leq q^R
     \text{ for every eigenvalue }\alpha\text{ of every }g_j\}.
\end{equation}
These sets are conjugation-invariant, open and closed. Every compact
subset of $L(F)$ is contained in some $B_R(L(F))$, by bounding the matrix
entries of $g_j$ and $g_j^{-1}$. Conversely, each $B_R(L(F))$ is bounded
modulo conjugation. Indeed, put a matrix in rational canonical form.
The coefficients of each invariant-factor polynomial are elementary
symmetric functions of a subset of its eigenvalues, hence are uniformly
bounded; the constant coefficients are uniformly bounded away from zero.
The corresponding companion matrices and their inverses therefore lie
in a fixed compact set. There are only finitely many possible block sizes. 

We also use the set $L(F)_{\mathrm{cm}}$ of elements compact modulo the
center. For $\GL(n,F)$ this means that all eigenvalues of $g$
have the same absolute value. Put
$$
 d(g)=(v_F(\det g_1),\ldots,v_F(\det g_s)),\qquad
 L(F)_{\mathrm{cm},k}=L(F)_{\mathrm{cm}}\cap d^{-1}(k),
 \quad k\in\mathbb Z^s.
$$
These sets are open and closed. Each $L(F)_{\mathrm{cm},k}$ is bounded
modulo conjugation: in its $j$th factor all eigenvalues have absolute
value $q^{-k_j/n_j}$.

For a conjugation-invariant open-and-closed subset $U\subset L(F)$,
multiplication by $\mathbf1_U$ induces a projector on
$\overline\mcH_L$. 
Denote the projectors for $B_R(L(F))$, $L(F)_{\mathrm{cm}}$, and
$L(F)_{\mathrm{cm},k}$ by $p_R$, $p_{\mathrm{cm}}$, and
$p_{\mathrm{cm},k}$, respectively.

Fix the standard Iwahori subgroup $I\subset\GL(n,F)$, and put
$$
 \mathfrak A=\{a\in M_n(\mathcal O_F):
                    a\bmod\varpi\text{ is upper triangular}\},
 \qquad I_m=1+\varpi^m\mathfrak A\quad(m\geq1).
$$
Thus the diagonal entries of an element of $I_m$ lie in
$1+\varpi^m\mathcal O_F$, its entries above the diagonal lie in
$\varpi^m\mathcal O_F$, and those below the diagonal lie in
$\varpi^{m+1}\mathcal O_F$. These are the integer-depth Moy--Prasad
subgroups at the barycenter of the standard alcove. For the group of points $L(F)$ of a product $L$,
take products of these subgroups. They form a basis of neighbourhoods of $1$.

For such a subgroup $K$, let $V_K$ be the image of $\mcH(L(F),K)$ in
$\overline\mcH_L$. All the projectors defined above preserve $V_K$;
we recall the precise finite-level statement in the proof below.
\begin{definition}\label{D:boundedsup} We say that a functional $\ell$ on $V_K$ has \emph{bounded support
modulo conjugation} if
\begin{equation}\label{E:bern-bounded-functional}
 \ell=\ell\circ p_R\quad\text{for some }R.
\end{equation}
\end{definition}

In particular, the restriction of any invariant distribution supported
in the conjugates of a compact set has this property.

\subsection{Induction and restriction} Choose representatives $\mathcal L(L)$
for the $L(F)$-conjugacy classes of Levi $F$-subgroups of $L$, including $L$ itself.
For $M\in\mathcal L(L)$, choose a parabolic subgroup $P=MU$ of $L$.
Then we have the restriction map
$$
 \bar r_M^L:\overline\mcH_L\longrightarrow
                 \overline\mcH_M
$$
adjoint to normalized parabolic induction. More precisely, for every
smooth admissible representation $\sigma$ of $M(F)$,
\begin{equation}\label{E:bern-restriction-trace}
 \tr\sigma(\bar r_M^L h)
 =\tr\bigl(i_{P(F)}^{L(F)}\sigma\bigr)(h).
\end{equation}
This map does not depend on $P$: the character on $L(F)$ of normalized induction from $M(F)$ is independent of the parabolic $F$-subgroup $P$ with Levi factor $M$.

In what follows we need a geometric description of this. Choose the Levi representatives to be standard, and choose each $P=MU$ to be standard.
Let $K_0=L(\mathcal O_F)$ be the standard maximal compact subgroup of $L(F)$, endowed with probability Haar measure $dk$. 
Normalize the Haar measure $du$ on $U(F)$ by
$$
 \operatorname{vol}_{U(F)}(K_0\cap U(F))
 =\frac{\operatorname{vol}_{L(F)}(K_0)}
        {\operatorname{vol}_{M(F)}(K_0\cap M(F))}.
$$
The class $\bar r_M^L h$ is represented by the normalized,
compact-averaged constant term
\begin{equation}\label{E:bern-constant-term}
 h^{P,K_0}(m)=|\delta_P(m)|^{1/2}
       \int_{K_0}\int_{U(F)}h(k^{-1}mu k)\,du\,dk,
 \qquad m\in M(F).
\end{equation}
Here $\delta_P$ is the character of $M$ given by the determinant of the adjoint action on the Lie algebra of $U$. We shall abbreviate
$\bar r_M^L$ to $\bar r_M$ when the ambient group is clear.

Thus an invariant distribution $\alpha_M$ on $M(F)$ induces the
invariant distribution $\alpha_M\circ\bar r_M$ on $L(F)$. In particular,
for $b\in\mcZ_M$ the induced distribution is
$$
 F_{M,b}(h)=b\bigl(\bar r_M^L(h)\bigr).
$$
Inversion commutes with restriction on cocenters:
$\bar r_M^L(h^\vee)=(\bar r_M^L h)^\vee$. This follows from
\eqref{E:bern-restriction-trace}, compatibility of normalized induction
with contragredients, and the density theorem for characters.

There is also a map of Bernstein centers
$\operatorname{res}_M:\mcZ_L\to\mcZ_M$, characterized by
\begin{equation}\label{E:bern-center-restriction}
 \bar r_M(b* h)=\operatorname{res}_M(b)*\bar r_M(h).
\end{equation}

\subsection{The discriminant and sufficiency}\label{SS:bern-statement}

Consider a Bernstein component of $G(F)=\GL(N,F)$
of the form
\begin{equation}\label{E:bern-component}
 \grs=
 \left[
 \prod_{a=1}^{r}\GL(d_a,F)^{m_a},
 \ \bigotimes_{a=1}^{r}\sigma_a^{\otimes m_a}
 \right],
 \qquad \sum_a d_am_a=N,
\end{equation}
where $\sigma_a$ is a unitary cuspidal representation of $\GL(d_a,F)$,
and these representations are pairwise
inequivalent up to unramified twist. Let
$$
 t_a=\#\{\chi\text{ unramified on }F^*:
          \sigma_a\otimes(\chi\circ\det)\simeq\sigma_a\}.
$$
Coordinates on the torus of distinct unramified twists are given by
$$
 z_{a,j}=\chi_{a,j}(\varpi)^{t_a},
 \qquad 1\leq j\leq m_a.
$$
Thus
$$
 Z_{\grs}
 =\bigotimes_a
 \mathbb C[z_{a,1}^{\pm1},\ldots,z_{a,m_a}^{\pm1}]^{S_{m_a}}.
$$
Define
\begin{equation}\label{E:bern-denominator}
 \Delta_{\grs}
 =\prod_a\prod_{\substack{1\leq i,j\leq m_a\\i\ne j}}
       \left(1-\frac{z_{a,i}}{z_{a,j}}\right).
\end{equation}
Changing $\sigma_a$ by an unramified twist multiplies all the coordinates
in its group by the same constant. Consequently
\eqref{E:bern-denominator} is independent of these choices. As usual we set the empty product to be 1.
These componentwise elements define
$$
 \Delta_{G}
 =\bigl(\Delta_{\grs}\bigr)_{\grs}
 \in\mcZ_G
 =\prod_{\grs} Z_{\grs}.
$$

For $L=\prod_j\GL(n_j)$, define
$\Delta_L\in\mcZ_L$ by taking the product of these discriminants
for its factors.

\begin{theorem}\label{T:denominator-main}
For every product $L$ of general linear groups,
$$
 \Delta_L*\mcD_L^c\subset\mcZ_L.
$$
Moreover, $\Delta_L$ is a non-zero-divisor in $\mcZ_L$.
Consequently Conjecture~\ref{conj:denominator} holds for $L$, and
$$
 (\Delta_L)\subset\mcJ_L,\qquad \mcD_L^c\subset\mcR_L.
$$
\end{theorem}

The formula for $\Delta_L$ is unchanged by contragredience: dualizing
the cuspidal representatives and inverting their twist coordinates
permutes its factors. In particular, $\Delta_L^\vee=\Delta_L$.
We prove the theorem using the following two statements.

\begin{theorem}\label{T:bern-main}
Let $L$ be a product of general linear groups.
\begin{enumerate}
\item\label{I:bern-decomposition}
\textup{(Decomposition over Levi subgroups.)}
Let $K\subset L(F)$ be one of the congruence subgroups specified above,
and let $\ell:V_K\to\mathbb C$ have bounded support modulo conjugation.
There are regular central elements
$b_M\in\mcZ_M$, each supported on finitely many Bernstein
components, such that
\begin{equation}\label{E:bern-Levi-decomposition}
 \ell(h)=\sum_{M\in\mathcal L(L)}
             b_M\bigl(\bar r_M^L(h)\bigr),
 \qquad h\in\mcH(L(F),K).
\end{equation}
The term $M=L$ is allowed. No discriminant occurs in this decomposition.

\item\label{I:bern-denominator}
Let $M$ be a Levi subgroup of $L$ and let
$b\in\mcZ_M$ have finite Bernstein support. For the induced
distribution
$$
 F_{M,b}(h)=b\bigl(\bar r_M^L(h)\bigr)
$$
there is a regular central element $A_{M,b}\in\mcZ_L$, with
finite Bernstein support, such that
\begin{equation}\label{E:bern-induced-denominator}
 F_{M,b}(\Delta_L* h)=A_{M,b}(h),
 \qquad h\in\mcH_L.
\end{equation}
\end{enumerate}
\end{theorem}

\begin{corollary}\label{C:bern-finite-level}
Let $L$ be a product of general linear groups and let
$K\subset L(F)$ be a congruence subgroup as above. For every functional
$\ell:V_K\to\mathbb C$ with bounded support modulo conjugation, the functional
$\Delta_L\ell$ is regular: there exists an element
$a_\ell\in\mcZ_L$, with finite Bernstein support, such that
\begin{equation}\label{E:bern-finite-level}
 (\Delta_L\ell)(h):=\ell(\Delta_L* h)
   =a_\ell(h),\qquad h\in\mcH(L(F),K).
\end{equation}
Thus multiplication by $\Delta_L$ makes every such functional come
from an element of the Bernstein center, viewed as a distribution.
\end{corollary}

\begin{proof}
Apply statement~(1) of Theorem~\ref{T:bern-main} to $\ell$, and statement~(2)
to each summand of \eqref{E:bern-Levi-decomposition}. Central operators
preserve $\mcH(L(F),K)$, so for every $h$ in this algebra,
$$
 \ell(\Delta_L* h)
 =\sum_{M\in\mathcal L(L)}A_{M,b_M}(h)
 =\left(\sum_{M\in\mathcal L(L)}A_{M,b_M}\right)(h).
$$
Thus $a_\ell=\sum_M A_{M,b_M}$ has the required properties.
\end{proof}

The proof of statement~(2) uses the type $A$ Plancherel formulas.
Statement~(1) uses separation of tempered representations by cuspidal
support, which is special to general linear groups. Neither argument
is asserted here for an arbitrary reductive group.
\subsection{Proof of Theorem~\ref{T:bern-main}\label{SS:bern-finite-level}(1)}

We use the following consequences of the finite-level Newton
decomposition; see \cite[Theorem 11 and Corollary 9]{He},
with the numbering of the published version:
\begin{enumerate}
\item The projections $p_R$, $p_{\mathrm{cm}}$, and
$p_{\mathrm{cm},k}$ preserve $V_K$. Their images, and the images of
their complementary projectors, have representatives in
$\mcH(L(F),K)$ supported in the indicated subsets.
\item The spaces $p_RV_K$ and $p_{\mathrm{cm},k}V_K$ are finite
dimensional, and
$$
 V_K=\bigcup_{R\geq0}p_RV_K,
 \qquad
 p_{\mathrm{cm}}V_K=\bigoplus_{k\in\mathbb Z^s}p_{\mathrm{cm},k}V_K.
$$
\end{enumerate}
For the first finite-dimensionality assertion one can equivalently
apply Howe's finiteness theorem to a compact set whose conjugates
cover $B_R(L(F))$; see \cite[Theorem 18]{He}. The assertion for
$L(F)_{\mathrm{cm},k}$ follows by its containment in some $B_R(L(F))$.

The statement below is crucial for our argument (and this is where we use heavily the fact that $L$ is a product of general linear groups). 
\begin{proposition}\label{P:bern-interpolation}
The functionals
$$
 V_K\longrightarrow\mathbb C,\qquad
 h\longmapsto b(h),\qquad b\in\mcZ_L,
$$
separate $V_K$; in other words if $h\in V_K$ is such that $b(h)=0$ for all $b\in \mcZ_L$, then $h=0$.
\end{proposition}
\begin{corollary}\label{C:fdim}
Let $W\subset V_K$ be a finite-dimensional subspace of $V_K$. Then the map $\mcZ_L\to W^*$ given by $b\mapsto(h\mapsto b(h))$ is surjective.
The element $b$ may be chosen to be supported on the finitely many Bernstein components visible at level $K$.
\end{corollary}
The corollary follows by finite-dimensional linear algebra. A central
element may be multiplied by the sum of the finitely many Bernstein
projectors visible at level $K$ without changing its restriction to $V_K$.
We now prove Proposition~\ref{P:bern-interpolation}.
\begin{proof}
The key observation is that an irreducible tempered representation of a
group $\GL(n,F)$ is determined by its cuspidal support. To see this,
work separately in each factor $\GL(n_j,F)$ of $L(F)$, and group the
cuspidal constituents of the cuspidal support by inertial equivalence.
For each class choose a unitary representative $\sigma_a$, and use the
coordinates $z=\chi(\varpi)^{t_a}$ introduced above, where $t_a$ is the
number of unramified self-twists of $\sigma_a$. The parameters in this
class form a multiset union of strings
\begin{equation}\label{E:bern-strings}
 uQ^{(\ell-1)/2},\ uQ^{(\ell-3)/2},\ldots,
 uQ^{-(\ell-1)/2},\qquad |u|=1,\qquad Q=q^{t_a}.
\end{equation}
The entries of largest absolute value determine the longest strings
and their unitary parameters. Removing these strings and repeating
recovers the complete multiset of strings. The classification of tempered
representations of $\GL(n,F)$, including irreducibility of
normalized induction of irreducible tempered representations, therefore
recovers the representation from its cuspidal support.

Let $X_K$ be the set of irreducible tempered representations of $L(F)$
with a nonzero $K$-fixed vector. Only finitely many Bernstein components
are visible at level $K$. Within each of them there are finitely many
tempered induction families, parametrized by compact tori modulo finite
groups. The function $\pi\mapsto\dim\pi^K$ is constant on each such
family, so $X_K$ is a finite union of these compact families.
The cuspidal-support map identifies $X_K$ with a compact subset of the
Bernstein variety. In particular, $X_K$ is Hausdorff and regular central
functions separate its points.

Their restrictions to $X_K$ are also closed under complex conjugation.
Indeed, every multiset of strings in \eqref{E:bern-strings} is invariant
under $z\mapsto\overline z^{-1}$. Thus, for a symmetric Laurent polynomial
$b$, conjugating its coefficients and replacing every coordinate by its
inverse gives a regular function whose restriction to $X_K$ is
$\overline b$.

The complex Stone--Weierstrass theorem says that a subalgebra of $C(X)$
on a compact Hausdorff space $X$ is uniformly dense if it contains the
constant functions, separates points, and is closed under complex
conjugation. These conditions hold for the restrictions of $\mcZ_L$
to $X_K$.

Now let $h\in\mcH(L(F),K)$ represent a class annihilated by all
$b\in\mcZ_L$, and put $f_h(\pi)=\tr\pi(h^\vee)$. This is a continuous function on
$X_K$, and it vanishes outside $X_K$. Plancherel gives
$$
 0=\int_{X_K} b(\pi)f_h(\pi)\,d\nu_L(\pi)\qquad(b\in\mcZ_L).
$$
The measure of $X_K$ is finite, since
$\int_{X_K}\dim\pi^K\,d\nu_L(\pi)=e_K(1)<\infty$.
Uniformly approximate $\overline{f_h}$ by regular central functions.
Passing to the limit yields
$$
 \int_{X_K}|f_h(\pi)|^2\,d\nu_L(\pi)=0.
$$
The Plancherel measure has full support on every tempered induction
family in $X_K$. Continuity therefore gives $f_h(\pi)=0$ at every point
of $X_K$, and hence on the entire tempered dual. To pass from tempered traces to all admissible traces, fix a tempered
induction datum and vary its unramified characters. The trace of $h^\vee$ is
a Laurent polynomial and vanishes on the unitary torus, so it vanishes
on the whole complex torus. The Langlands classification expresses every
irreducible character as a linear combination of such standard-module
characters. Spectral density, in the form of
\cite[Theorem 2.8]{CH}, now implies that the class of $h^\vee$,
and hence of $h$, in $\overline\mcH_L$ is zero. This form applies
over non-Archimedean local fields of either characteristic.
For the split groups used here, see also
\cite[Theorem 0]{K-cusp} in characteristic zero and
\cite{K-close} for passage to positive characteristic.
\end{proof}

We prove Theorem~\ref{T:bern-main}(1) by induction on the semisimple rank of $L$. The following statement will be needed in this induction argument.

\begin{lemma}\label{L:bern-finite-support}
For every $b\in\mcZ_L$, its restriction to $V_K$ has bounded support
modulo conjugation.
\end{lemma}

\begin{proof}
The function $b_K=b* e_K$ is compactly supported. Choose $R$ with
$\operatorname{supp}(b_K)\subset B_R(L(F))$. A class in
$(1-p_R)V_K$ has a representative $h\in\mcH(L(F),K)$ supported outside
$B_R(L(F))$. For this representative,
$$
 b(h)=\int_{L(F)} b_K(g)h(g)\,dg=0.
$$
Thus $b|_{V_K}=(b|_{V_K})\circ p_R$.
\end{proof}

We now use the following cocenter theorem of Ciubotaru--He:
\begin{equation}\label{E:bern-elliptic-support}
 \bigcap_{M\subsetneq L}\ker\bar r_M
 \ \subset\ p_{\mathrm{cm}}\overline\mcH_L,
\end{equation}
where $M$ ranges over $\mathcal L(L)\setminus\{L\}$.
In words, a cocenter class whose restrictions to all proper Levi
subgroups vanish is represented by functions supported on elements
compact modulo the center; see \cite[Theorem 4.4]{CH}.

There is also the elementary compatibility
\begin{equation}\label{E:bern-bounded-restriction}
 \bar r_Mp_R^L=p_R^M\bar r_M.
\end{equation}
Indeed, the compact-averaged constant-term formula defining
$\bar r_M$ evaluates $h$ on conjugates of $mu$, with $m\in M(F)$ and
$u\in U(F)$, where $U$ is the unipotent radical of $P$. A block upper-triangular matrix $mu$
has precisely the eigenvalues of its diagonal blocks $m$, so
$mu\in B_R(L(F))$ if and only if $m\in B_R(M(F))$.

\begin{lemma}\label{L:bern-descent}
Let $\ell$ be a functional on $V_K$ with bounded support modulo
conjugation, and suppose that $\ell$ vanishes on $p_{\mathrm{cm}}V_K$.
Then, for suitable congruence subgroups $K_M\subset M(F)$, there are
functionals $\ell_M$ on $V_{K_M}$, each with bounded support modulo
conjugation, such that
\begin{equation}\label{E:bern-descent}
 \ell=\sum_{M\in\mathcal L(L),\ M\ne L}\ell_M\circ\bar r_M
 \quad\text{on }V_K.
\end{equation}
\end{lemma}

\begin{proof}
Choose $R$ such that $\ell=\ell\circ p_R$. Choose $K_M$ sufficiently
small that $\bar r_M(V_K)\subset V_{K_M}$ for each $M$. A uniform
level exists by the compact-averaged constant-term formula: first
pass to a principal congruence subgroup normal in a maximal compact
subgroup and then to a smaller subgroup of the prescribed form in $M(F)$.

Consider the map of finite-dimensional spaces
$$
 p_R^LV_K\longrightarrow
       \bigoplus_{M\subsetneq L}p_R^MV_{K_M},
 \qquad h\longmapsto(\bar r_Mh)_M.
$$
By \eqref{E:bern-elliptic-support}, its kernel is contained in
$p_{\mathrm{cm}}V_K$, so $\ell$ annihilates that kernel. Thus $\ell$
factors through the image of this map. Extend the resulting functional
to the displayed direct sum and write it as a sum of functionals
$\lambda_M$ on the individual factors. Set
$\ell_M=\lambda_M\circ p_R^M$ on $V_{K_M}$.
Equation \eqref{E:bern-bounded-restriction} gives, for every $h\in V_K$,
$$
 \ell(h)=\ell(p_R^Lh)
 =\sum_M\lambda_M(\bar r_Mp_R^Lh)
 =\sum_M\ell_M(\bar r_Mh).
$$
Each $\ell_M$ satisfies \eqref{E:bern-bounded-functional}, as required.
\end{proof}

\subsection{End of the proof of Theorem~\ref{T:bern-main}(1)}
We use induction on the semisimple rank of $L$.
If $L$ is a torus, the bounded-support condition means that the
functional on $\mcH(L(F),K)$ is supported on finitely many cosets
in $L(F)/K$. It is therefore represented, using the fixed Haar
measure, by a compactly supported $K$-invariant function.
As $L(F)$ is abelian, this defines a central distribution; its
Bernstein support is contained in the finitely many components
visible at level $K$. This proves the base case.

Suppose now that the assertion holds for proper Levi subgroups of $L$.
We first match $\ell$ on all elements compact modulo the center.
If $\ell=\ell\circ p_R$, its restriction to
$p_{\mathrm{cm},k}V_K$ can be nonzero only for the finitely many
multidegrees satisfying $|k_j|\leq n_jR$ for every $j$.
For each such $k$, Corollary~\ref{C:fdim} applied to $W=p_{\mathrm{cm},k}V_K$ implies that there exists some 
$b_k\in\mcZ_L$ with finite Bernstein support such that
$$
 b_k(h)=\ell(h),\qquad
 h\in p_{\mathrm{cm},k}V_K.
$$
Retain only the homogeneous component of $b_k$ of determinant
multidegree $k$. More explicitly, for $\mathbf c\in(\mathbb C^\times)^s$
put
$$
 \nu_{\mathbf c}(g)=\prod_j c_j^{v_F(\det g_j)},\qquad
 (T_{\mathbf c}b)(h)=b(\nu_{\mathbf c}h).
$$
Its spectral function is $(T_{\mathbf c}b)(\pi)
=b(\pi\otimes\nu_{\mathbf c})$. For a regular central element
$b$ with finite Bernstein support the algebraic twisting action
has a finite Laurent expansion
$$
 T_{\mathbf c}b=\sum_{k\in\mathbb Z^s}\mathbf c^k b^{[k]},
 \qquad b^{[k]}\in\mcZ_L.
$$
If $h$ is supported in determinant multidegree $j$, then
$(T_{\mathbf c}b)(h)=\mathbf c^j b(h)$.
Comparison of coefficients gives
$b^{[k]}(h)=\delta_{k,j}b(h)$.
Thus replacing $b_k$ by $b_k^{[k]}$ preserves the displayed equality,
retains finite Bernstein support, and makes the distribution vanish
on every other determinant multidegree.

Set $b_L=\sum_k b_k$. Then
$$
 b_L(h)=\ell(h),\qquad h\in p_{\mathrm{cm}}V_K.
$$
Consequently
$\eta(h)=\ell(h)-b_L(h)$ vanishes on
$p_{\mathrm{cm}}V_K$. It still has bounded support modulo conjugation,
by Lemma~\ref{L:bern-finite-support}. Lemma~\ref{L:bern-descent} gives
$$
 \eta=\sum_{M\in\mathcal L(L),\ M\ne L}
             \ell_M\circ\bar r_M^L,
$$
where $\ell_M$ has bounded support on a suitable $V_{K_M}$.
Apply the induction hypothesis to every $\ell_M$. Each is a sum of
distributions induced from Bernstein center elements of Levi subgroups $D$ of $M$. Transitivity,
$$
 \bar r_D^M\circ\bar r_M^L=\bar r_D^L,
$$
turns these into distributions induced from Bernstein center elements of Levi subgroups of $L$.
Transporting conjugate Levi subgroups to the chosen representatives
and adding their central elements gives
\eqref{E:bern-Levi-decomposition}. There are only finitely many terms,
and every central element obtained has finite Bernstein support which finishes the proof.
\subsection{Proof of Theorem~\ref{T:bern-main}\label{SS:bern-transfer}(2)}
We use an explicit calculation with the Plancherel formula.
First observe that the quotient
\begin{equation}\label{E:bern-cross-D}
 C_M=\frac{\operatorname{res}_M(\Delta_L)}{\Delta_M}
\end{equation}
is a regular element of $\mcZ_M$. Since the discriminants are invariant
under inversion and restriction commutes with inversion, $C_M^\vee=C_M$.
We now claim
\begin{lemma}\label{L:bern-transfer}
For $a\in\mcZ_M$ supported on finitely many Bernstein components,
there is $\mathcal T_M(a)\in\mcZ_L$, also supported on finitely many
components, such that
\begin{equation}\label{E:bern-transfer}
 (a*C_M)\bigl(\bar r_M(h)\bigr)
 =\mathcal T_M(a)(h),
 \qquad h\in\mcH_L.
\end{equation}
\end{lemma}

\begin{proof}
Fix a component $\grt$ of $M(F)$ mapping to a component $\grs$ of $L(F)$.
For a product $L$, the following indices and cross pairs are taken
separately within each of its general linear factors.
For each inertial class represented by $\sigma_a$, write $m_{a,b}$ for the number of its
occurrences assigned to the $b$th $M$-factor. Put
$$
 Q_a=q^{t_a},\qquad
 W_\grs=\prod_a S_{m_a},\qquad
 W_\grt=\prod_{a,b}S_{m_{a,b}},
$$
and define the regular function
\begin{equation}\label{E:bern-cross-Q}
 P_{\grs/\grt}
 =\prod_a
  \prod_{\substack{i,j\text{ assigned to}\\
                     \text{different }M\text{-factors}}}
       \left(1-Q_a^{-1}\frac{z_{a,i}}{z_{a,j}}\right).
\end{equation}
Both orientations of every cross pair occur in this product.

For each type $A$ affine Hecke algebra use the canonical trace
$\tau(T_w)=\delta_{w,1}$ in its standard basis, and probability
Haar measure on its compact unramified-parameter torus.
The Plancherel measure of the component $\grs$ is a constant
$\kappa_{L,\grs}>0$ times the product of these normalized type $A$
affine-Hecke Plancherel measures with parameters $Q_a$.
The same holds for $(M,\grt)$, with constant
$\kappa_{M,\grt}>0$. These comparison constants incorporate the
fixed Haar measures on the groups. This includes all tempered induction families,
by \cite[Theorems 5.7 and 6.12]{AP}.

Set
$$
 c_\grt
 =\frac{\kappa_{M,\grt}}{\kappa_{L,\grs}}
   \prod_a Q_a^{\binom{m_a}{2}-\sum_b\binom{m_{a,b}}{2}}.
$$
The contribution of $a_\grt$ to $\mathcal T_M(a)|_\grs$ is
\begin{equation}\label{E:bern-transfer-polynomial}
 c_\grt\sum_{w\in W_\grs/W_\grt}
        w\bigl(a_\grt P_{\grs/\grt}\bigr).
\end{equation}
It is regular and $W_\grs$-invariant.

We verify the distribution identity by applying Plancherel to $h^\vee$
and using $\bar r_M(h^\vee)=(\bar r_Mh)^\vee$. On each tempered family, the parameters
are unions of the strings \eqref{E:bern-strings}. If a placement of the
coordinates splits a string between two $M$-factors, some adjacent
pair in that string lies in different factors. Their ratio is $Q_a$
in one orientation, so the corresponding factor in
\eqref{E:bern-cross-Q} vanishes. Thus only placements of entire strings
in $M$-factors contribute.

For a surviving placement, pass to a common parameter cover of the
tempered family on the two sides. The relative type $A$ Plancherel
factor is
$$
 \left(\prod_a
 Q_a^{-\binom{m_a}{2}+\sum_b\binom{m_{a,b}}{2}}\right)
 \frac{C_M}{P_{\grs/\grt}}.
$$
The discrete-series factors are the same on this common cover.
We spell out the finite-group multiplicities, making the following
count separately in each general linear factor of $L$.
For each $a$, let $n_{a,\ell}$ be the number of strings of length
$\ell$ in the chosen tempered family, and put
$$
 A=\prod_{a,\ell}n_{a,\ell}!.
$$
Fix an assignment type of entire strings to the $M$-factors, and let
$n_{a,b,\ell}$ be the number of length-$\ell$ strings of class $a$
assigned to factor $b$. Thus
$$
 \sum_b n_{a,b,\ell}=n_{a,\ell},
 \qquad
 \sum_\ell \ell\,n_{a,b,\ell}=m_{a,b}.
$$
Put $B=\prod_{a,b,\ell}n_{a,b,\ell}!$. On the ordered torus of
string parameters, the Weyl factors in the Plancherel measures of
the $L$-family and this $M$-family are respectively $1/A$ and $1/B$.
Exactly $A/B$ cosets in $W_\grs/W_\grt$ give placements of this
assignment type. Permuting strings of equal length within each
class $a$ carries any such placement to any other. Although the
corresponding values of $a_\grt$ need not agree pointwise, the
integrals of the corresponding summands in
\eqref{E:bern-transfer-polynomial} agree after this relabelling:
the ordered Haar measure, the $L$-Plancherel density, and the
induced $L$-character are invariant under these permutations.
Consequently the coset sum converts the factor $1/A$ into $1/B$,
as required.

The power of $Q_a$ is independent of the string partition.
Indeed, if the strings in class $a$ have lengths $\ell_{a,i}$,
the difference between the exponents in the two parabolic
Plancherel factors is
$$
 \left(\binom{m_a}{2}
       -\sum_i\binom{\ell_{a,i}}{2}\right)
 -
 \left(\sum_b\binom{m_{a,b}}{2}
       -\sum_i\binom{\ell_{a,i}}{2}\right)
 =
 \binom{m_a}{2}-\sum_b\binom{m_{a,b}}{2}.
$$
This is precisely the exponent occurring in $c_\grt$.
The relative density formula used here is the parabolic
factorization of Plancherel measure
\cite[Proposition 3.27 and formula (4.96)]{Opdam}.
Together with the componentwise constants
$\kappa_{M,\grt}/\kappa_{L,\grs}$, it shows that multiplication by
\eqref{E:bern-transfer-polynomial} gives exactly the contribution
of $(a*C_M)(\bar r_Mh)$ from the chosen assignment type.
For a product $L$, take the product of these multiplicity and
density comparisons. Summing over assignment types gives the
identity on the chosen tempered family. These comparisons are
made where the string parameters are in general position;
the omitted subset has measure zero in each family. Summing
over all tempered families and the finitely many relevant
Bernstein components proves \eqref{E:bern-transfer}.
\end{proof}

\begin{proof}[Proof of Theorem~\ref{T:bern-main}(2)]
For $b\in\mcZ_M$ as in the statement, set
$$
 A_{M,b}=\mathcal T_M(b\Delta_M)\in\mcZ_L.
$$
This element is regular and has finite Bernstein support. By
\eqref{E:bern-center-restriction}, \eqref{E:bern-cross-D}, and
Lemma~\ref{L:bern-transfer}, we obtain
\begin{align*}
 F_{M,b}(\Delta_L* h)
 &=b\bigl(\operatorname{res}_M(\Delta_L)*\bar r_Mh\bigr)\\
 &=(b*\Delta_M*C_M)(\bar r_Mh)\\
 &=\mathcal T_M(b\Delta_M)(h).
\end{align*}
Here we used $\Delta_M^\vee=\Delta_M$ and $C_M^\vee=C_M$ to pass
from convolution on test functions to convolution on distributions.
This is \eqref{E:bern-induced-denominator}. When $M=L$, the assertion
simply reads $A_{L,b}=b\Delta_L$. For a proper Levi subgroup, it is
the Plancherel comparison in Lemma~\ref{L:bern-transfer} that supplies
the regular numerator.
\end{proof}

\subsection{Passage from finite levels to invariant distributions}
\begin{proof}[Proof of Theorem~\ref{T:denominator-main}]
Let $\alpha\in\mcD_L^c$. Its restriction to every $V_K$ has bounded
support in the sense of \eqref{E:bern-bounded-functional}.
By Corollary~\ref{C:bern-finite-level} there is $a_K\in\mcZ_L$ such that
\begin{equation}\label{eq:localtrace}
 \alpha(\Delta_L*h)=a_K(h),
 \qquad h\in\mcH(L(F),K).
\end{equation}
Put $\beta=\Delta_L*\alpha$. Since $\Delta_L^\vee=\Delta_L$, the
definition of the central action gives
$$
 \beta(h)=\alpha(\Delta_L*h)=a_K(h),
 \qquad h\in\mcH(L(F),K).
$$
The convolution $v_K=\beta*e_K$ is a locally constant $K$-bi-invariant
function; compact support is not needed to define this convolution
with a test function. For $h\in\mcH(L(F),K)$, \eqref{eq:localtrace} becomes
$$
 \int_{L(F)}v_K(g)h(g)\,dg
 =\int_{L(F)}(a_K*e_K)(g)h(g)\,dg.
$$
Characteristic functions of individual $K$-double cosets show that
$v_K=a_K*e_K$. Thus $\beta*e_K$ has compact support.
The chosen congruence subgroups form a neighbourhood basis of the
identity. Every test function is bi-invariant under one of them, and
therefore $\beta*h=(\beta*e_K)*h$ is compactly supported for every
$h\in\mcH_L$. The distribution $\beta$ is invariant and essentially
compact, hence belongs to $\mcZ_L$.

Finally, each $Z_{\grs}$ is a domain and the displayed product defining
$\Delta_{\grs}$ is not identically zero. Thus $\Delta_L$ is a
non-zero-divisor in the product ring $\mcZ_L$. This proves all the
assertions.
\end{proof}

\begin{corollary}\label{cor:orbital}
Let $y\in\GL(N,F)$ be regular semisimple and let $\mu_y$ be an
invariant measure on its conjugacy class. Then
$$
 \Delta_G*\mu_y\in\mcZ_G.
$$
For every Bernstein component $\grs$ there is a regular element
$A_{y,\grs}\in Z_{\grs}$ such that
$$
 \mu_y*(\Delta_{\grs}*f)=A_{y,\grs}*f,
 \qquad f\in(\mcH_G)_{\grs}.
$$
Moreover, $\Delta_G$ is nonzero at every irreducible representation
of $G(F)$ that is square-integrable modulo the center.
\end{corollary}
\begin{proof}
The orbital measure belongs to $\mcD_G^c$, so the compactness assertion
follows from the theorem. For the final assertion, such a representation
is a generalized Steinberg representation $\St(\sigma,m)$, with
$N=dm$. For a unitary $\sigma$ its cuspidal support has twist
coordinates
$$
 z_j=uQ^{(m+1-2j)/2},\qquad Q=q^{t_\sigma}>1,\quad |u|=1.
$$
Consequently
$$
 \Delta_G(\St(\sigma,m))=\prod_{i\ne j}(1-Q^{j-i})\ne0.
$$
Unramified twists multiply all these coordinates by the same scalar
and leave the product unchanged.
\end{proof}

\subsection{Necessity of the discriminant}

\begin{proof}[Completion of the proof of Theorem~\ref{T:denominator-exact}]
We give the proof for $G=\GL(N)$; for a product, the same construction
is performed in each factor. Fix a Bernstein component $\grs$ with
cuspidal support Levi
$$
 M=\prod_a\GL(d_a)^{m_a},\qquad
 \sigma=\bigotimes_a\sigma_a^{\otimes m_a}.
$$
Let $\grt=[M(F),\sigma]_M$ be the corresponding cuspidal component
of $M(F)$, and let $e_{\grt}\in\mcZ_M$ be its Bernstein projector.
Define the invariant distribution induced from this projector by
$$
 \alpha_{\grs}(h)=e_{\grt}\bigl(\bar r_M^G h\bigr).
$$

The projector $e_{\grt}$ is supported on compact elements of $M(F)$;
see \cite[Theorem 1.3 and Section 1.6]{BK}. In the constant-term
formula for $\bar r_M^G$, an element $mu$ has the same eigenvalues
as its block-diagonal part $m$. Hence a test function vanishing on
$B_0(G(F))$ has constant term vanishing on $B_0(M(F))$.
It follows that $\alpha_{\grs}$ is supported in $B_0(G(F))$. This set is bounded modulo conjugation,
as proved after \eqref{E:bern-bounded-set}. Thus
$\alpha_{\grs}\in\mcD_G^c$.

On $\grt$ one has $\Delta_M=1$ and $W_{\grt}=1$. Put
$$
 P_{\grs}=\prod_a\prod_{i\ne j}
       \left(1-Q_a^{-1}\frac{z_{a,i}}{z_{a,j}}\right),
 \qquad Q_a=q^{t_a}>1.
$$
Formula~\eqref{E:bern-transfer-polynomial}, with $a=e_{\grt}$,
therefore gives
\begin{equation}\label{E:sharp-trace}
 \Delta_G*\alpha_{\grs}=c_{\grs}P_{\grs},\qquad c_{\grs}\ne0,
 \quad\text{as invariant distributions,}
\end{equation}
where the central element on the right is extended by zero off
$\grs$. The constant includes $|W_{\grs}|$, since $P_{\grs}$ is
$W_{\grs}$-invariant.

Now let $z\in\mcJ_G$. By definition,
$b=z*\alpha_{\grs}\in\mcZ_G$. Commutativity of the central action
and \eqref{E:sharp-trace} imply, on the component $\grs$,
\begin{equation}\label{E:sharp-divisibility}
 \Delta_{\grs}b_{\grs}
       =c_{\grs}z_{\grs}P_{\grs}.
\end{equation}
The ring $Z_{\grs}$ is a Laurent polynomial ring in elementary
symmetric functions, and in particular is a unique factorization
domain. The discriminant $\Delta_{\grs}$ is squarefree in this ring;
its prime factors are the reduced collision divisors, one for each
$a$ with $m_a\geq2$ (see Proposition~\ref{prop:reduced-discriminant}).
At the generic point of any such divisor, precisely one pair of
coordinates coincides. Its two factors in $P_{\grs}$ are
$1-Q_a^{-1}$, which is nonzero, and the other factors are generically
nonzero. Thus $P_{\grs}$ is relatively prime to $\Delta_{\grs}$.
Equation~\eqref{E:sharp-divisibility} gives
$\Delta_{\grs}\mid z_{\grs}$.

This holds on every component. The quotients
$(z_{\grs}/\Delta_{\grs})_{\grs}$ define an element of
$\mcZ_G=\prod_{\grs}Z_{\grs}$, so $z\in\Delta_G\mcZ_G$.
The reverse inclusion is Theorem~\ref{T:denominator-main}.
For products, take $M$, $e_{\grt}$, and $P_{\grs}$ as the corresponding
products; the same coprimality and divisibility argument applies.
\end{proof}

\section{The endoscopic intersection for general linear groups}
\label{sec:intersection-typeA}
\subsection{The reduced collision locus in type A}
For a Bernstein component as in \eqref{E:bern-component}, put
$$
 T_{\grs}=\prod_a(\mathbb C^\times)^{m_a},\qquad
 W_{\grs}=\prod_a S_{m_a},\qquad
 X_{\grs}=T_{\grs}/W_{\grs}.
$$
Let $D_{\grs}\subset X_{\grs}$ be the reduced image of the union of
the hyperplanes $z_{a,i}=z_{a,j}$, $i\ne j$.
\begin{proposition}\label{prop:reduced-discriminant}
The ideal of $D_{\grs}$ is the principal ideal $(\Delta_{\grs})$.
In particular, the common denominator of
Theorem~\ref{T:denominator-main} defines the reduced collision divisor
on the quotient, although its pullback has order two on every
collision hyperplane of the ordered torus.
\end{proposition}
\begin{proof}
Work first in one symmetric group factor. If an invariant Laurent
polynomial $f$ vanishes on $z_i=z_j$, it is divisible by $z_i-z_j$.
Write $f=(z_i-z_j)g$. Invariance under the transposition $(ij)$ makes
$g$ anti-invariant. Its restriction to $z_i=z_j$ is therefore zero,
so $f$ is divisible by $(z_i-z_j)^2$. Applying this to every pair and
using unique factorization in the Laurent polynomial ring gives
$$
 \prod_{i<j}(z_i-z_j)^2\mid f.
$$
Up to an invertible symmetric Laurent monomial and a nonzero constant,
this product is $\prod_{i\ne j}(1-z_i/z_j)$. The quotient by this
product is again invariant. Repeating the argument for each $a$ proves
the assertion.
\end{proof}

Under local Langlands for $\GL(N)$, the points on this divisor have
repeated isomorphic irreducible Weil constituents in their semisimple
parameters. If two such constituents are placed in different blocks of
an endoscopic Levi, their off-diagonal intertwining spaces increase
the centralizer dimension. Thus the divisor has the centralizer-jump
interpretation motivating Section~\ref{sec:endoscopy}.

\subsection{The Fitting ideals for endoscopic Levis}
For $\GL(N)$, ordinary untwisted endoscopic transfers are the standard
block-sum transfers from products of general linear groups, up to
equivalence and changes of source coordinates. Indeed, in the
untwisted case the Weil action fixes the semisimple endoscopic element;
its centralizer is a product of general linear groups, and the action
on that product is inner. This assertion would not apply to endoscopy
with a nontrivial twisting character.

Let $H=\prod_b\GL(n_b)$ and fix a source component mapping to
$X_{\grs}$. Assign each irreducible Weil constituent indexed by
$(a,i)$ to its $H$-block $b(a,i)$. This describes the source quotient as
$T_{\grs}/W_{\grt}$, where
$W_{\grt}=\prod_{a,b}S_{m_{a,b}}$. Its map to $X_{\grs}$ is finite.
For fixed $\grs$ only finitely many block partitions and assignments
occur. The same description applies on each bounded parameter piece.
Under the local Langlands correspondence, these quotients are the
corresponding Bernstein components. The map $\Psi_G$ is this
identification, by \cite[Theorem IX.7.4]{FS} and the characteristic-zero
description of the invariant parameter rings in
\cite[Theorem 1.7 and Corollary 1.8]{DHKM}. The latter is an
isomorphism of schemes in characteristic zero, so this identification
is not limited to closed geometric points.

Define on the source component
\begin{equation}\label{eq:cross-collision}
 C_{\grs/\grt}
  =\prod_a\prod_{\substack{i,j\\b(a,i)\ne b(a,j)}}
          \left(1-\frac{z_{a,i}}{z_{a,j}}\right).
\end{equation}
Both orientations of every pair occur. With the compatible
coordinates of Section~\ref{SS:bern-transfer}, this is the restriction
of the analytic cross factor $C_M$ in \eqref{E:bern-cross-D} to
the source component. The internal-block
discriminant $\Delta_{\grt}$ satisfies
\begin{equation}\label{eq:geometric-factorization}
 f^*\Delta_{\grs}=\Delta_{\grt}C_{\grs/\grt}.
\end{equation}

\begin{lemma}\label{lem:cross-fitting}
On the parameter scheme over this source component, the pullback of
$C_{\grs/\grt}$ belongs to the zeroth Fitting ideal of the universal
coinvariant module for
$V=\mathfrak{gl}_N/\bigoplus_b\mathfrak{gl}_{n_b}$.
Consequently $f^*\Delta_{\grs}$ belongs to the contracted ideal
$J_{\mathfrak e,U}$ on this component.
\end{lemma}
\begin{proof}
We work in characteristic zero. Locally on fixed-inertial-type
charts, write the universal parameter in Weil--Deligne form
$(r,\mathcal N)$, with $r$ having finite inertia image and
$r(\mathrm{Fr})=\varphi(\mathrm{Fr})$, where $\varphi$ is the
universal Weil parameter. These are scheme-theoretic charts:
see \cite[Propositions 2.6--2.7, Section 3, and (4.2)--(4.5)]{DHKM}.
In particular, the argument below applies to the universal
family, not only to geometric points.
After a faithfully flat change of charts, we may fix the finite
inertia action. Let $R$ be the chart ring and put
$$
 W=(R\otimes V^\vee)_{r(I_F)}.
$$
Taking coinvariants for a finite group in characteristic zero is
a direct-summand operation, so $W$ is locally free. Frobenius
acts on $W$.

Let $\mathcal N_V$ be the nilpotent operator induced by monodromy on
$R\otimes V^\vee$. A suitable power of a tame inertia generator
acts as $\exp(c\mathcal N_V)$ for some nonzero scalar $c$. Since
$$
 \exp(c\mathcal N_V)-1=\mathcal N_Vu(\mathcal N_V),\qquad u(0)=c\ne0,
$$
the operator $u(\mathcal N_V)$ is invertible. Thus its coinvariant relation
imposes the $\mathcal N_V$-relation. The original inertia relations then
impose the finite-inertia relations. Consequently the full
coinvariant module $M$ is obtained by imposing the monodromy
relations on
$$
 \operatorname{coker}(1-r(\mathrm{Fr}):W\longrightarrow W).
$$
In particular, $M$ is a quotient of this cokernel.
Zeroth Fitting ideals increase under taking quotients. Therefore
\begin{equation}\label{eq:fitting-determinant}
 \det(1-r(\mathrm{Fr})\mid W)\in\Fitt_0^R(M).
\end{equation}

We compute this determinant using the irreducible constituents of
the semisimplification of $r$. Finite-inertia coinvariants are exact,
and determinants are multiplicative along invariant filtrations,
so this calculation also applies when Frobenius is not semisimple.
Two irreducible Weil representations have a nonzero inertia
intertwiner precisely when they differ by an unramified twist.
For the class indexed by $a$, the restriction of its representative
to inertia has $t_a$ pairwise nonisomorphic constituents, cyclically
permuted by Frobenius. For two twists in this class, Frobenius acts
cyclically on the corresponding $t_a$-dimensional inertia-intertwiner
space, with cycle product $z_{a,i}/z_{a,j}$. Its contribution to
$\det(1-r(\mathrm{Fr}))$ is therefore
$1-z_{a,i}/z_{a,j}$, with the orientation reversed on the dual.
The module $V$ consists exactly of the cross-block Hom spaces,
and both orientations occur. Thus the determinant in
\eqref{eq:fitting-determinant} is $C_{\grs/\grt}$, up to a unit.
This calculation is an identity in the chart ring: it is the
characteristic-polynomial identity for the finite inertia types,
not merely an equality on a generic open subset.

Fitting ideals commute with base change, and ideal membership
descends under faithfully flat maps. This proves the first assertion
on the parameter scheme. Multiplying by $\Delta_{\grt}$ and using
\eqref{eq:geometric-factorization} proves the second assertion after
contraction to the invariant ring. In particular, monodromy directions
do not invalidate the containment.
\end{proof}

\subsection{Computation of the intersection}
\begin{theorem}\label{T:intersection-typeA}
For $G=\GL(N)$ and every Bernstein component $\grs$, the endoscopic
intersection, under the spectral-to-Bernstein identification, is
$$
 \Icap_{\grs}=(\Delta_{\grs}).
$$
Globally,
$$
 {\mcZ_G\cdot\Psi_G(\Icap_G)
            =\Delta_G\mcZ_G=\mcJ_G.}
$$
Thus Conjecture~\ref{conj:geometric} holds for $\GL(N)$.
\end{theorem}
\begin{proof}
Lemma~\ref{lem:cross-fitting} shows that $\Delta_{\grs}$ belongs
to every endoscopic ideal on this component. Hence
$$
 (\Delta_{\grs})\subset\Icap_{\grs}.
$$
For the reverse inclusion, use the finest support Levi
$$
 H_{\grs}=\prod_a\GL(d_a)^{m_a}
$$
and the source component in which each factor has one irreducible
Weil constituent. Its coarse parameter space is the ordered torus
$T_{\grs}$. One may take an additional finite \'etale cover
$\widetilde T_{\grs}\to T_{\grs}$ to choose the individual
unramified twists before quotienting by their self-twists. On this
cover there is an actual family of block-diagonal Weil parameters
with these irreducible constituents and with $\mathcal N=0$.

For this family there are no monodromy relations and no repeated
constituents inside any one block. The universal coinvariant module
is exactly the cokernel of $1-r(\mathrm{Fr})$ on the finite-inertia
coinvariants. Its zeroth Fitting ideal on $\widetilde T_{\grs}$
is consequently
$$
 (\Delta_{\grs})\OO(\widetilde T_{\grs}).
$$
This is a calculation of the full ideal, including at multiple
collisions. If $a\in\Icap_{\grs}$, its pullback to this family
must lie in that ideal, since Fitting ideals commute with base
change. The map
$\widetilde T_{\grs}\to X_{\grs}$ is faithfully flat:
the self-twist power maps are finite \'etale, and the symmetric
Laurent-polynomial extension for each $S_{m_a}$ is finite free.
Contraction of a principal extended ideal therefore gives
$$
 (\Delta_{\grs})\OO(\widetilde T_{\grs})
       \cap\OO(X_{\grs})=(\Delta_{\grs}).
$$
It follows that $a\in(\Delta_{\grs})$. If the support Levi is
$G$ itself, then $\Delta_{\grs}=1$, and the first inclusion already
gives the assertion without using a proper transfer.

The formulas use only the multiset of irreducible constituents,
are independent of coordinate representatives, and agree upon
changing ramification bounds. They therefore compute the global
ideal defined by \eqref{eq:global-intersection} as the product of
these component ideals. The spectral-to-Bernstein identification
and Theorem~\ref{T:denominator-exact} give the global equality.
Finally $\Delta_G$ is a non-zero-divisor, as in
Theorem~\ref{T:denominator-main}.
\end{proof}

\section{The conjectures for \texorpdfstring{$\mathrm{SL}(2)$}{SL(2)}}\label{sec:sl2}
Put $G=\mathrm{SL}(2)$ and $\widetilde G=\GL(2)$. We prove both
Conjectures~\ref{conj:denominator} and \ref{conj:geometric} for $G$.
The additional phenomenon, compared with $\GL(2)$, is the splitting
of a principal series at a nontrivial quadratic character. Differences
of the two constituent characters force additional zeros in the
denominator. The nonsplit endoscopic tori account for these zeros.

\subsection{The statement and its coordinates}
Identify the diagonal torus $T(F)$ with $F^\times$ by
$a\mapsto\operatorname{diag}(a,a^{-1})$, and write $I(\chi)$ for
normalized induction of a smooth character $\chi$ of $T(F)$ to
$G(F)$. A principal series Bernstein component is specified by
$\vartheta=\chi|_{\mathcal O_F^\times}$, up to inversion.
Choose an extension $\chi_0$ of $\vartheta$ and write
$\chi=\chi_0\nu_z$, where $\nu_z$ is unramified and
$\nu_z(\varpi)=z$. If $\vartheta^2=1$, choose
$\chi_0(\varpi)=1$, so that $\chi_0^2=1$ and inversion acts by
$z\mapsto z^{-1}$. The component rings are
$$
 Z_{\grs}=\begin{cases}
   \mathbb C[z^{\pm1}],&\vartheta^2\ne1,\\
   \mathbb C[z^{\pm1}]^{z\mapsto z^{-1}}=\mathbb C[t],
       \quad t=z+z^{-1},&\vartheta^2=1.
 \end{cases}
$$
A cuspidal Bernstein component is a point, since $G(F)$ has no
nontrivial unramified characters. Define $d\in\mcZ_G$ componentwise by
\begin{equation}\label{eq:sl2-denominator}
 d_{\grs}=\begin{cases}
  1,&\grs\text{ is cuspidal},\\
  1,&\grs\text{ is principal series and }\vartheta^2\ne1,\\
  t^2-4,&\grs\text{ is principal series and }\vartheta^2=1.
 \end{cases}
\end{equation}
Changing the quadratic extension $\chi_0$ replaces $t$ by $-t$,
so the last formula is independent of that choice.

\begin{theorem}\label{T:sl2}
For $G=\mathrm{SL}(2)$,
$$
 {\mcJ_G=d\mcZ_G
          =\mcZ_G\cdot\Psi_G(\Icap_G).}
$$
The element $d$ is a non-zero-divisor. In particular,
$d*\mcD_G^c\subset\mcZ_G$ and $\mcD_G^c\subset\mcR_G$.
\end{theorem}

We use the standard restriction theory from $\GL(2,F)$ to
$\mathrm{SL}(2,F)$: restriction of an irreducible representation
is a finite sum of irreducibles, its constituents form a single
$\GL(2,F)$-orbit, and every irreducible representation of
$\mathrm{SL}(2,F)$ occurs this way. Cuspidality is preserved
under restriction. See \cite[Section 2 and Theorem 3.1]{Tadic}.
In the principal series, $I(\chi)$ is irreducible except when
$\chi=|\cdot|^{\pm1}$ or $\chi$ is nontrivial quadratic.
In the latter case
$$
 I(\chi)=\pi_\chi^+\oplus\pi_\chi^-,
$$
with the two inequivalent constituents interchanged by
$\GL(2,F)$-conjugation. At $|\cdot|^{\pm1}$ the constituents are
the trivial and Steinberg representations, each fixed by this
conjugation action. These assertions also follow directly from the
restriction theory: an irreducible principal series of $\GL(2,F)$
has a nontrivial determinant self-twist precisely when its character
ratio is nontrivial quadratic; see \cite[Sections 2--3]{Tadic}.

\subsection{Descent of the \texorpdfstring{$\GL(2)$}{GL(2)} denominator}
Let $j:G(F)\hookrightarrow\widetilde G(F)$ be the inclusion. Define
$d_0\in\mcZ_G$ to be $2-t=(1-z)(1-z^{-1})$ on the unramified
principal series component and $1$ on every other component.
It is invariant under $\widetilde G(F)$-conjugation.

\begin{lemma}\label{L:sl2-descent}
If $\alpha\in\mcD_G^c$ is invariant under $\widetilde G(F)$,
then $d_0*\alpha\in\mcZ_G$.
\end{lemma}
\begin{proof}
For a distribution $\gamma$ on $G(F)$ let
$j_*\gamma(f)=\gamma(f|_{G(F)})$. An outer-invariant central
distribution $b\in\mcZ_G$ gives a central distribution
$j_*b\in\mcZ_{\widetilde G}$. To check essential compactness,
use the homeomorphism
$$
 F^\times\times G(F)\longrightarrow\widetilde G(F),
 \qquad (a,h)\longmapsto s(a)h,
 \quad s(a)=\operatorname{diag}(a,1).
$$
Every compactly supported locally constant function in these
coordinates is a finite sum of functions $\phi(a)f_0(h)$.
Convolution by $j_*b$ sends such a function to
$\phi(a)(b*f_0)(h)$, because $b$ is invariant under conjugation
by $s(a)$. This is compactly supported. The same observation shows
conversely that, for an outer-invariant distribution $\gamma$,
essential compactness of $j_*\gamma$ implies essential compactness
of $\gamma$: take $\phi(1)\ne0$ and restrict the convolution to
the fiber $a=1$.

The central element $j_*d_0$ is precisely
$\Delta_{\widetilde G}$ of \eqref{E:bern-denominator}.
Indeed, on a principal series of $\widetilde G(F)$ with inducing
characters $\chi_1,\chi_2$, restriction to $G(F)$ has character
ratio $\chi_1\chi_2^{-1}$. The $\GL(2)$ discriminant is
$(1-z)(1-z^{-1})$ when this ratio is unramified, and is $1$
otherwise. On cuspidal representations both elements act by
$1$. The same scalar comparison holds on every constituent at
reducibility points, and hence identifies the central elements.

Now $j_*\alpha$ is invariant under $\widetilde G(F)$ and has bounded
support modulo conjugation there. Theorem~\ref{T:denominator-main}
gives
$$
 j_*(d_0*\alpha)=\Delta_{\widetilde G}*j_*\alpha
                       \in\mcZ_{\widetilde G}.
$$
The converse compactness observation proves the lemma.
\end{proof}

\subsection{The part not invariant under \texorpdfstring{$\GL(2,F)$}{GL(2,F)}}
The conjugation action on invariant distributions and on the cocenter
of $G(F)$ factors through
$$
 C=\widetilde G(F)/(Z_{\widetilde G}(F)G(F))
                 \simeq F^\times/(F^\times)^2.
$$
This is a compact group. It is finite if $\operatorname{char}F\ne2$;
in characteristic $2$ we use normalized Haar averaging on $C$.
The action on every test function modulo commutators is smooth, so
averaging is well-defined. Write $\operatorname{Av}$ for this
projector, on either cocenters or invariant distributions.
It preserves bounded support modulo conjugation. Indeed, $C$ has
representatives $s(a)$ with
$a\in\mathcal O_F^\times\cup\varpi\mathcal O_F^\times$, a compact
set, and conjugating a fixed compact support set by these
representatives again gives a compact set.

For nontrivial quadratic $\chi$, put
$$
 \eta_\chi(h)
   =\tr\pi_\chi^+(h^\vee)-\tr\pi_\chi^-(h^\vee).
$$
Let $x_\chi$ be its common central character in the principal series
Bernstein component. Then
$u*\eta_\chi=u(x_\chi)\eta_\chi$ for $u\in\mcZ_G$.

\begin{lemma}\label{L:sl2-noninvariant}
For any invariant distribution $\beta$ with
$\operatorname{Av}(\beta)=0$, its projection to a principal series
component is a linear combination of the distributions $\eta_\chi$
in that component. Its projections to the cuspidal components
together form an element of $\mcZ_G$.
\end{lemma}
\begin{proof}
Apply the trace description to $\ell_\beta(h)=\beta(h^\vee)$.
On a principal series component consider the image of
$1-\operatorname{Av}$ in the cocenter. The character of any
$\widetilde G(F)$-fixed irreducible representation vanishes on this
image. By the preceding classification, only the two constituents
at each nontrivial quadratic character remain, and their traces on
this image are opposite. There are at most two such characters in
any one component. More explicitly, let
$\overline\mcH_{G,\grs}$ denote the image of the Hecke algebra of
the component $\grs$ in $\overline\mcH_G$. Spectral density
\cite[Theorem 2.8]{CH} gives an injective map
$$
 \begin{aligned}
 (1-\operatorname{Av})\overline\mcH_{G,\grs}
 &\longrightarrow
 \bigoplus_{\substack{\chi^2=1,\ \chi\ne1\\\chi\in\grs}}\mathbb C,\\
 [h]&\longmapsto
 \bigl(\tr\pi_\chi^+(h)-\tr\pi_\chi^-(h)\bigr)_\chi.
 \end{aligned}
$$
Here $\chi\in\grs$ means that the cuspidal support of $I(\chi)$
belongs to $\grs$. Indeed, if these differences vanish, then all
irreducible traces vanish on the given class, so density makes
that class zero. Every functional on the image of this map
extends to its finite-dimensional target. Dualizing therefore
gives the assertion for $\ell_\beta$, and inversion gives the
assertion for $\beta$.

Each cuspidal component contains a single irreducible
representation $\pi$. Density makes its cocenter one-dimensional.
By Plancherel, the Bernstein projector $e_\pi$, regarded as a distribution,
satisfies
$$
 e_\pi(h)=d_\pi\,\tr\pi(h^\vee),
$$
where $d_\pi$ is the nonzero formal degree of $\pi$.
Thus every invariant distribution on this component is a scalar multiple
of $e_\pi$ and belongs to $\mcZ_G$.
Arbitrary scalars on the cuspidal components define an element
of the product Bernstein center; only finitely many components are
visible on any given compact-open level. This proves the last
assertion.
\end{proof}

Define $r\in\mcZ_G$ by
$$
 r_{\grs}=\begin{cases}
 t+2,&\vartheta=1,\\
 t^2-4,&\vartheta\ne1,\ \vartheta^2=1,\\
 1,&\text{on all other components}.
 \end{cases}
$$
It vanishes at every $x_\chi$ with $\chi$ nontrivial quadratic.
Consequently Lemma~\ref{L:sl2-noninvariant} gives
$r*(\alpha-\operatorname{Av}(\alpha))\in\mcZ_G$ for every
invariant distribution $\alpha$. If $\alpha\in\mcD_G^c$,
Lemma~\ref{L:sl2-descent} also gives
$d_0*\operatorname{Av}(\alpha)\in\mcZ_G$. Combining these,
$$
 d_0r*\alpha\in\mcZ_G.
$$
The element $d_0r$ differs from $d$ only by the unit $-1$ on the
unramified principal series component. Thus
\begin{equation}\label{eq:sl2-sufficiency}
 d\mcZ_G\subset\mcJ_G.
\end{equation}
Every component of $d$ is nonzero, so this already proves
Conjecture~\ref{conj:denominator} for $\mathrm{SL}(2)$.

\subsection{Necessity of the factors}
First, $\eta_\chi\in\mcD_G^c$. To see this, on the split regular
semisimple set every $\GL(2,F)$-conjugacy class is one
$\mathrm{SL}(2,F)$-conjugacy class: the determinant map on the
split centralizer in $\GL(2,F)$ is surjective. The two constituent
characters therefore agree there. In particular, their difference
vanishes outside the set
$$
 B_0(G(F))=\{g\in G(F):\text{both eigenvalues have absolute value }1\}.
$$
Indeed, every element outside this set is split regular semisimple.
This set is bounded modulo $G(F)$-conjugation. Rational canonical
form first proves boundedness modulo $\GL(2,F)$-conjugation, and
the compact representatives of $C$ used above convert this into
boundedness modulo $G(F)$-conjugation.

The nonzero distribution $\eta_\chi$ is not central. If it were
represented by $b\in Z_{\grs}$, then
$(a-a(x_\chi))b=0$ for every $a\in Z_{\grs}$. Since
$Z_{\grs}=\mathbb C[t]$ is a domain of positive dimension, this
would force $b=0$. Hence, for $u\in\mcJ_G$,
$$
 u(x_\chi)=0\qquad(\chi\ne1,\ \chi^2=1).
$$
This forces $t+2$ on the unramified component and $t^2-4$ on every
ramified quadratic component.

It remains to force $t-2$ on the unramified component. Let $e_T$ be
the Bernstein projector to the unramified characters of $T(F)$ and
put
$$
 \alpha_T(h)=e_T\bigl(\bar r_T^G h\bigr).
$$
The projector $e_T$ is supported on $T(\mathcal O_F)$.
The constant-term formula therefore shows that
$\alpha_T\in\mcD_G^c$, supported in $B_0(G(F))$.
The rank-one Plancherel formula gives
\begin{equation}\label{eq:sl2-plancherel}
 d_0*\alpha_T=cP,\qquad
 P(t)=(1-q^{-1}z)(1-q^{-1}z^{-1})
      =1+q^{-2}-q^{-1}t,
 \quad c\ne0,
\end{equation}
where $P$ is extended by zero off the unramified component.
For completeness, the continuous Plancherel density on its unitary
torus is a nonzero constant times
$(1-z)(1-z^{-1})/P(t)$; see
\cite[Proposition 3.12, Theorems 3.14 and 4.5(2)(iii), and Sections 4.7--4.8]{KM},
specialized to the unramified case $q_1=q_2=q$.
Multiplication by $P$ leaves the torus trace weighted by $d_0$.
The same formulas show that the only discrete-series contribution
in this component is Steinberg, with $t=q+q^{-1}$, where $P=0$.
This proves \eqref{eq:sl2-plancherel}; the isolated reducible unitary
principal series contributes no Plancherel atom.

For $u\in\mcJ_G$ write $b=u*\alpha_T\in\mcZ_G$. On the
unramified component, \eqref{eq:sl2-plancherel} gives
$$
 (2-t)b=c\,uP(t).
$$
Since $P(2)=(1-q^{-1})^2\ne0$, one has $(t-2)\mid u$.
Together with the quadratic-character factors, this proves
$\mcJ_G\subset d\mcZ_G$. Combining with
\eqref{eq:sl2-sufficiency} proves the first equality in
Theorem~\ref{T:sl2}.

\subsection{The endoscopic ideal}
The dual group is $\widehat G=\mathrm{PGL}(2)$.
Besides the datum with dual group $\widehat G$, which imposes no
condition, the connected endoscopic dual groups are tori.
Their Weil action is either trivial or inversion. Thus the proper
endoscopic groups are the split torus and the norm-one tori of
quadratic extensions $E/F$.

For the split torus, the quotient
$\Lie(\widehat G)/\Lie(\widehat T)$ is the sum of the two root
lines. On the family $\chi=\chi_0\nu_z$ their Weil characters are
$\chi$ and $\chi^{-1}$. If $\vartheta\ne1$, finite-inertia
coinvariants vanish, so this datum imposes no condition. If
$\vartheta=1$, the zeroth Fitting ideal on the source torus is
$$
 ((1-z)(1-z^{-1})),
$$
whose contraction to $\mathbb C[t]$ is $(t-2)$.

Now let $H$ be the norm-one torus of $E/F$, and let
$\kappa_E$ be the corresponding quadratic character of $W_F$,
or of $F^\times$ by local class field theory. Its $L$-embedding
has image in the normalizer of $\widehat T$, with nontrivial
Weyl action specified by $\kappa_E$. Its source coarse parameter
components are reduced points. Indeed, $H(F)$ is compact, and
\cite[Proposition IX.6.4]{FS} identifies its spectral center with
$$
 \varprojlim_K\Lambda[H(F)/K].
$$
For compact open $K\subset H(F)$, the quotient $H(F)/K$ is finite
abelian, so its group algebra over $\Lambda$ is a product of
copies of $\Lambda$. Together with the characteristic-zero
coarse-quotient description in \cite[Theorem 1.7]{DHKM}, this
shows that the bounded coarse parameter pieces are finite and
reduced.

The Lie algebra of the centralizer
inside $\widehat H$ is zero, since Weil acts on it by inversion.
The centralizer inside $\widehat G$ has positive-dimensional Lie
algebra precisely when the transferred parameter has image the
order-two group generated by a Weyl reflection. Indeed, the
centralizer of a reflection has connected component a torus;
requiring the whole parameter to centralize that torus leaves
only this order-two image. Up to conjugacy this is the parameter
$$
 w\longmapsto\operatorname{diag}(\kappa_E(w),1)
                       \quad\text{in }\mathrm{PGL}(2).
$$
We apply the construction of Section~\ref{sec:endoscopy} on one
such source component. Its coordinate ring is $S=\Lambda$.
Write $R$ for the coordinate ring of the parameter scheme above
it and $M$ for the universal coinvariant module. The ideal used
on the source quotient is
$$
 J=\Fitt_0^R(M)\cap S.
$$
At the exceptional source point, the fiber of $M$ is nonzero,
so $\Fitt_0^R(M)$ is proper. Its contraction to the field $S$
is therefore zero. Hence the annihilator of the pushforward of
$S/J=\Lambda$ on the target is exactly the maximal ideal of
$x_{\kappa_E}$.

At every other source point, the closed-orbit representative has
zero coinvariant space. The invariant-support argument of
Section~\ref{sec:algebraic} shows that the support has empty
image on this component. Thus the contracted ideal is the unit
ideal. In particular, the exceptional contribution is a reduced
point, with multiplicity one.

It follows that the endoscopic intersection on the unramified
principal series component is
$$
 (t-2)\cap(t+2)=(t^2-4).
$$
On a ramified quadratic component the two quadratic extensions of
$\vartheta$ give the ideals $(t-2)$ and $(t+2)$, again with
intersection $(t^2-4)$. On all other principal series components
the ideal is the unit ideal. At a cuspidal parameter the
centralizer in $\mathrm{PGL}(2)$ has zero Lie algebra, and hence
there is no exceptional locus; the ideal is again the unit ideal.
These are calculations of the Fitting and annihilator ideals, not
only of their supports.

To compare with $\mcZ_G$, the Fargues--Scholze parameter of a
constituent of a restricted $\GL(2,F)$-representation is the
projectivization of its semisimple parameter, by
\cite[Theorems IX.6.1 and IX.7.4]{FS}. On principal series
components this gives exactly the coordinates used above.
Several cuspidal Bernstein points can have the same spectral
parameter; this causes no change to the extended ideal, since the
spectral ideal is the unit ideal at every such parameter.
The componentwise functions just calculated define a spectral
element whose image is $d$, and every element of the spectral
ideal pulls back to a multiple of $d$. Thus
$$
 \mcZ_G\cdot\Psi_G(\Icap_G)=d\mcZ_G.
$$
The calculations are independent of ramification bounds and of
the representatives of the endoscopic data. They prove the
second equality of Theorem~\ref{T:sl2} and complete the proof of
Conjecture~\ref{conj:geometric} for $\mathrm{SL}(2)$.


\clearpage
\begin{thebibliography}{99}
\bibitem{Tadic}
M. Tadi\'c, \emph{Notes on representations of non-Archimedean
$\mathrm{SL}(n)$}, Pacific Journal of Mathematics \textbf{152}
(1992), no.~2, 375--396.
\href{https://msp.org/pjm/1992/152-2/pjm-v152-n2-p11-s.pdf}{Publisher's text}.

\bibitem{BK}
A. Braverman and D. Kazhdan, with an appendix by R. Bezrukavnikov,
\emph{Bernstein components via Bernstein center},
\href{https://arxiv.org/abs/1512.08637}{arXiv:1512.08637},
version 3 (2018).
\bibitem{FS}
L. Fargues and P. Scholze, \emph{Geometrization of the local Langlands correspondence}, Ast\'erisque \textbf{466} (2026), Soci\'et\'e Math\'ematique de France. See also arXiv:2102.13459, version 4 (2024).
\bibitem{H}
D. Hansen, \emph{Excursion operators and the stable Bernstein center}, Forum of Mathematics, Pi \textbf{14} (2026), e10, 1--14. DOI: 10.1017/fmp.2026.10028.
\bibitem{Varma}
S. Varma, \emph{Some comments on the stable Bernstein center},
preprint, revision of July 3, 2024.
\href{https://mathweb.tifr.res.in/~sandeepv/stable_center_classical_groups_revision_July_3_2024.pdf}{Author's manuscript}.
\bibitem{Cotner}
S. Cotner, \emph{Morphisms of character varieties}, International Mathematics Research Notices \textbf{2024}, no.~16, 11540 ff.; arXiv:2304.10135, version 3 (2024), Theorem~1.1.
\bibitem{BGM}
M. Bate, H. Geranios, and B. Martin, \emph{Orbit closures and invariants}, arXiv:1604.00924, version 2 (2018), Section~2.2, Lemma~2.3.
\bibitem{DHKM}
J.-F. Dat, D. Helm, R. Kurinczuk, and G. Moss, \emph{Moduli of Langlands parameters}, Journal of the European Mathematical Society \textbf{27} (2025), no.~5, 1827--1927, DOI: 10.4171/JEMS/1599; arXiv:2009.06708, version 3 (2024).
See especially Theorem~1.7, Corollary~1.8, Proposition~2.7,
and Sections~3--4.
\bibitem{AP}
A.-M. Aubert and R. Plymen,
\emph{Plancherel measure for $\GL(n,F)$ and $\GL(m,D)$:
explicit formulas and Bernstein decomposition},
Journal of Number Theory \textbf{112} (2005), 26--66.
\href{https://arxiv.org/abs/math/0302169}{arXiv:math/0302169}.

\bibitem{CH}
D. Ciubotaru and X. He,
\emph{Cocenters of $p$-adic groups, III: elliptic and rigid cocenters},
Peking Mathematical Journal \textbf{4} (2021), 159--186.
\href{https://arxiv.org/abs/1703.00378}{arXiv:1703.00378}.

\bibitem{He}
X. He,
\emph{Cocenters of $p$-adic groups, I: Newton decomposition},
Forum of Mathematics, Pi \textbf{6} (2018), e2.
\href{https://arxiv.org/abs/1610.04791}{arXiv:1610.04791}.

\bibitem{K-cusp}
D.~Kazhdan, {\em Cuspidal geometry of \(p\)-adic groups}, Journal d'Analyse Math\'ematique {\bf 47} (1986), 1--36.

\bibitem{K-close}
D.~Kazhdan, \emph{Representations of groups over close local fields},
Journal d'Analyse Math\'ematique \textbf{47} (1986), 175--179.

\bibitem{KM}
P. Kutzko and L. Morris,
\emph{Explicit Plancherel Theorems for $H(q_1,q_2)$ and
$\mathrm{SL}_2(F)$},
Pure and Applied Mathematics Quarterly \textbf{5} (2009),
no.~1, 435--467.
\href{https://archive.intlpress.com/site/pub/files/_fulltext/journals/pamq/2009/0005/0001/PAMQ-2009-0005-0001-a013.pdf}{Publisher's text}.

\bibitem{Opdam}
E. M. Opdam,
\emph{On the spectral decomposition of affine Hecke algebras},
Journal of the Institute of Mathematics of Jussieu \textbf{3} (2004),
531--648.
\href{https://arxiv.org/abs/math/0101007}{arXiv:math/0101007}.
\end{thebibliography}
\end{document}